\documentclass[english,12pt,a4paper,reqno]{amsart}
\usepackage{amssymb,amsmath,amstext,amsthm,amsfonts}
\usepackage[ansinew]{inputenc}
\usepackage{a4wide}

\newcommand{\um}{{\bf 1}}

\newcommand{\leb}{\operatorname{Leb}}

\newcommand{\dist}{\operatorname{dist}}
\newcommand{\diam}{\operatorname{diam}}

\newcommand{\ve}{\varepsilon}
\newcommand{\vf}{\varphi}

\newcommand{\norm}[1]{\left\Vert#1\right\Vert}
\newcommand{\mdlo}[1]{\left| #1\right|}
\newcommand{\ch}[1]{\left\{ #1\right\}}

\begin{document}

\newcommand{\mcup}{\mbox{$\bigcup$}}
\newcommand{\mcap}{\mbox{$\bigcap$}}

\def \RR {{\mathbb R}}
\def \ZZ {{\mathbb Z}}
\def \NN {{\mathbb N}}
\def \PP {{\mathbb P}}
\def \TT {{\mathbb T}}
\def \II {{\mathbb I}}
\def \JJ {{\mathbb J}}

\def \vare {\varepsilon }

 \def \cf {\mathcal{F}}
 \def \cm {\mathcal{M}}
 \def \cn {\mathcal{N}}
 \def \cq {\mathcal{Q}}
 \def \cp {\mathcal{P}}
 \def \cb {\mathcal{B}}
 \def \cc {\mathcal{C}}
 \def \cs {\mathcal{S}}
 \def \bc {\mathcal{B}}
 \def \hc {\mathcal{F}}
 \def \ca {\mathcal{A}}

\newcommand{\dem}{\begin{proof}}
\newcommand{\cqd}{\end{proof}}

\newcommand{\qand}{\quad\text{and}\quad}

\newtheorem{theorem}{Theorem}
\newtheorem{corollary}{Corollary}

\newtheorem*{Maintheorem}{Main Theorem}

\newtheorem{maintheorem}{Theorem}
\renewcommand{\themaintheorem}{\Alph{maintheorem}}
\newcommand{\cmt}{\begin{maintheorem}}
\newcommand{\fmt}{\end{maintheorem}}

\newtheorem{maincorollary}[maintheorem]{Corollary}
\renewcommand{\themaintheorem}{\Alph{maintheorem}}
\newcommand{\cmc}{\begin{maincorollary}}
\newcommand{\fmc}{\end{maincorollary}}

\newtheorem{T}{Theorem}[section]
\newcommand{\ct}{\begin{T}}
\newcommand{\ft}{\end{T}}

\newtheorem{Corollary}[T]{Corollary}
\newcommand{\cco}{\begin{Corollary}}
\newcommand{\fco}{\end{Corollary}}

\newtheorem{Proposition}[T]{Proposition}
\newcommand{\cpr}{\begin{Proposition}}
\newcommand{\fpr}{\end{Proposition}}

\newtheorem{Lemma}[T]{Lemma}
\newcommand{\cle}{\begin{Lemma}}
\newcommand{\fle}{\end{Lemma}}

\newtheorem{Sublemma}[T]{Sublemma}
\newcommand{\csle}{\begin{Lemma}}
\newcommand{\fsle}{\end{Lemma}}

\theoremstyle{remark}

\newtheorem{Remark}[T]{Remark}
\newcommand{\cre}{\begin{Remark}}
\newcommand{\fre}{\end{Remark}}

\newtheorem{Definition}{Definition}
\newcommand{\cd}{\begin{Definition}}
\newcommand{\fd}{\end{Definition}}

\title[Recurrence times and rates of mixing]{Recurrence times and rates of mixing\\ for invertible dynamical systems}

\author{José F. Alves}
\address{José F. Alves \\Departamento de Matemática Pura, Faculdade de Ciências da Universidade do Porto\\
Rua do Campo Alegre 687, 4169-007 Porto, Portugal}
\email{jfalves@fc.up.pt}
\urladdr{http://www.fc.up.pt/cmup/home/jfalves}

\author{Vilton  Pinheiro}
\address{Vilton  Pinheiro \\Departamento de Matemática, Universidade Federal da Bahia\\
Av. Ademar de Barros s/n, 40170-110 Salvador, Brazil.}
\email{viltonj@ufba.br}

\date{\today}


\thanks{Work carried out at the Federal University of
Bahia, University of Porto and IMPA. JFA was partially supported by
FCT through CMUP and POCI/MAT/61237/2004. VP was partially
supported by PADCT/CNPq and POCI/MAT/61237/2004}


\begin{abstract}

We consider invertible discrete-time dynamical systems having a hyperbolic product structure in some region of the phase space with infinitely many branches and variable recurrence time. We show that the decay of correlations of the SRB measure associated to that hyperbolic structure is related to the tail of the recurrence times. We also give sufficient conditions for the validity of the Central Limit Theorem. This generalizes previous results by Benedicks and Young.
\end{abstract}

\maketitle

\setcounter{tocdepth}{1}

\tableofcontents


\section{Introduction}

One of the most powerful ways of describing the dynamical features of
chaotic dynamical systems  is through invariant probability measures.
A map \( f \) is said to
be \emph{mixing} with respect to an invariant probability  measure $\mu$ if
\[
|\mu(f^{-n}(A) \cap B) - \mu(A)\mu(B)|\to 0,\quad\text{when
$n\to\infty$},
\]
for any measurable sets \( A, B \).
Standard counterexamples show that in general there is no specific
rate at which this convergence to 0 occurs. However, defining the
\emph{correlation function} of {\em observables} $\varphi,\psi\colon M\to\RR$,
\[
\mathcal C_{n}(\varphi, \psi;\mu) = \left|\int (\varphi \circ f^{n}) \psi d\mu - \int \varphi d\mu \int \psi
d\mu\right|,
\]
it is sometimes possible  to  obtain specific rates of decay, which depend only on the map~\( f \) (up to a multiplicative constant which is allowed to
depend on \( \varphi, \psi \)), provided the observables
\( \varphi, \psi \) have sufficient  regularity.
Notice that choosing these observables to be characteristic
functions this gives
exactly the definition
of mixing.
Still in this direction, the {\em Central Limit
Theorem} states that the probability of a given deviation of
the average values of an observable along an orbit from the
spatial  average is essentially given by a Normal Distribution. 

Since the work of Sinai, Ruelle and Bowen \cite{S,R,B} it is known that  uniformly hyperbolic diffeomorphisms (Axiom A, Anosov) possess {\em SRB (or physical) measures} with exponential decay of correlations and satisfying  the Central Limit
Theorem. By physical  measure we mean an invariant probability measure such that for a large set (positive volume) of initial states the asymptotic time average (with respect to a continuous observable) coincides with the spatial average of that observable (with respect to the measure). A key ingredient in the proofs of Sinai, Ruelle and Bowen are Markov partitions, which permit to codify the dynamics and from its codification to deduce the main statistical features of the dynamical system.

In the context of non-uniformly hyperbolic diffeomorphisms, Benedicks and Young introduced in \cite{BY2} some kind of 
 structures with Markov flavor in certain regions of the phase space  with infinitely many branches
and variable return times.  This structures enabled them to obtain exponential decay of correlations and deduce the Central Limit Theorem
for Hénon maps.
Further developments by Young in \cite{Y1} lead to a
joint treatment of some non-uniformly hyperbolic diffeomorphisms, including Hénon maps, billiards with convex scatterers and Axiom~A attractors. This kind of approach has also been successfully implemented  by
Young in \cite{Y2} for studying the rates of mixing of non-invertible systems with some non-uniformly expanding behavior.

The frameworks developed by Young in  \cite{Y1} and \cite{Y2} are certainly among the most
powerful tools for studying the statistical properties of non-uniformly hyperbolic dynamical systems.
In both approaches, there is an explicit relation between the tail of the recurrence times to the hyperbolic structure and the decay of correlations, at least
for some specific rates. However, the results in both
papers do not depict reasonably the whole scenario. On the one hand, the model
in \cite{Y1} can only be applied to systems whose decay of
correlations is exponential. On the other hand, the model in
\cite{Y2}, in spite of being suitable for other decay rates,  is
specific to non-invertible systems. Let us mention that such a simple diffeomorphism  as the solenoid with intermittency that we present in Section~\ref{se.application} does not fit the model in \cite{Y1}; see Remark~\ref{re.hypti}.

The present work essentially aims at being a step farther in the construction of a theory on the  statistical features of non-uniformly hyperbolic diffeomorphisms. We believe that hyperbolic structures with sub-exponential tail of recurrence times can play an important role in obtaining the rates of mixing for the diffeomorphisms introduced by Viana in \cite{V}. Such hyperbolic structures can possibly be useful also in the study of some classes of billiards and Poincaré return maps for flows, for which the tails of recurrence frequently decay at sub-exponential rates.

 \subsection{Hyperbolic structures}\label{se.hypstructures}
Let $f\colon M\to M$ be defined on a finite
dimensional Riemannian manifold $M$, and let $\leb$ denote a
normalized volume form on the Borel sets of $M$ that we call it {\em
Lebesgue measure.} Given a submanifold $\gamma\subset M$   we use
$\leb_\gamma$ to denote the measure on $\gamma$ induced by the
restriction of the Riemannian structure to~$\gamma$.

An embedded disk $\gamma\subset M$ is called an {\em unstable
manifold} if $\dist(f^{-n}(x),f^{-n}(y))\to0$ exponentially fast as
$n\to\infty$ for every $x,y\in\gamma$. Similarly, $\gamma$ is called
a {\em stable manifold} if $\dist(f^{n}(x),f^{n}(y))\to0$
exponentially fast as $n\to\infty$ for every $x,y\in\gamma$.

\cd Let  $\text{Emb}^1(D^u,M)$ be the space of $C^1$ embeddings from
$D^u$ into $M$. We say that $\Gamma^u=\{\gamma^u\}$ is a
\emph{continuous family of $C^1$ unstable manifolds} if there is a
compact set~$K^s$, a unit disk $D^u$ of some $\RR^n$, and a map
$\Phi^u\colon K^s\times D^u\to M$ such that
\begin{itemize}
\item[i)] $\gamma^u=\Phi^u(\{x\}\times D^u)$ is an unstable
manifold;
\item[ii)] $\Phi^u$ maps $K^s\times D^u$ homeomorphically onto its
image; \item[iii)] $x\mapsto \Phi^u\vert(\{x\}\times D^u)$ defines a
continuous map from $K^s$ into $\text{Emb}^1(D^u,M)$.
\end{itemize}
Continuous families of $C^1$ stable manifolds are defined similarly.

\fd


\cd We say that $\Lambda\subset M$ has a \emph{hyperbolic product
structure} if there exist a continuous family of unstable
manifolds $\Gamma^u=\{\gamma^u\}$ and a continuous family of
stable manifolds $\Gamma^s=\{\gamma^s\}$ such that
\begin{itemize}
    \item[i)] $\Lambda=(\cup \gamma^u)\cap(\cup\gamma^s)$;
    \item[ii)] $\dim \gamma^u+\dim \gamma^s=\dim M$;
    \item[iii)] each  $\gamma^s$ meets each $\gamma^u$  in exactly one point;
    \item[iv)] stable and unstable manifolds are transversal with angles bounded
    away from~0.
\end{itemize}
\fd

Let $\Lambda\subset M$ have a
hyperbolic product structure, whose defining families are $\Gamma^s$
and $\Gamma^u$. A subset $\Lambda_0\subset \Lambda$ is called an
{\em $s$-subset} if $\Lambda_0$ also has a hyperbolic product
structure and its defining families $\Gamma_0^s$ and $\Gamma_0^u$
can be chosen with $\Gamma_0^s\subset\Gamma^s$ and
$\Gamma_0^u=\Gamma^u$; {\em $u$-subsets} are defined analogously.
Given $x\in\Lambda$, let $\gamma^{*}(x)$ denote the element of
$\Gamma^{*}$ containing $x$, for $*=s,u$. For each $n\ge 1$ let
$(f^n)^u$ denote the restriction of the map $f^n$ to
$\gamma^u$-disks, and let $\det D(f^n)^u$ be the Jacobian of
$D(f^n)^u$.
We require that the  hyperbolic product structure~$\Lambda$ satisfies several properties:

 {\em
\begin{enumerate}
    \item[\bf
    (P$_1$)] {\em Markov:} there are pairwise disjoint $s$-subsets $\Lambda_1,\Lambda_2,\dots\subset\Lambda$ such
    that
    \begin{enumerate}
 \item $\leb_{\gamma}\big((\Lambda\setminus\cup\Lambda_i)\cap\gamma\big)=0$ on each $\gamma\in\Gamma^u$;
 \item for each $i\in\NN$ there is $R_i\in\NN$ such that $f^{R_i}(\Lambda_i)$ is $u$-subset,
         and for all $x\in \Lambda_i$
         $$
         f^{R_i}(\gamma^s(x))\subset \gamma^s(f^{R_i}(x))\qand
         f^{R_i}(\gamma^u(x))\supset \gamma^u(f^{R_i}(x)).
         $$

    \end{enumerate}
\end{enumerate}
}

In the  statements of  the  remaining properties about the hyperbolic structure we assume that $C>0$ and
$0<\beta<1$ are constants  which only depend on $f$ and $\Lambda$.
\smallskip

{\em
 \begin{enumerate}
\item[\bf(P$_2$)] {\em Contraction on  stable leaves}:
  $\dist(f^n(y),f^n(x))\le C\beta^n $, $\forall y\in\gamma^s(x)\,\forall n\ge
 1$.
 \end{enumerate}
 }
 
 \smallskip
 
In spite of the uniform contraction in the stable direction, this condition is not too restrictive in systems having regions where the contraction fails to be uniform, since we are allowed to remove points in the unstable leaves, provided a subset with positive measure in those leaves remains at the end. This has been carried out in \cite{BY2} for Hénon maps.


Next we introduce a return time function $R\colon \Lambda\to\NN$ and a
return map $f^R\colon \Lambda\to\Lambda$,  defined for each
$i\in\NN$ as $$
 R\vert_{\Lambda_i}=R_i\qand
f^R\vert_{\Lambda_i}=f^{R_i}\vert_{\Lambda_i}. $$
         We consider
the {\em separation time} $s(x,y)$ for $x,y\in\Lambda$ as
$$
s(x,y)=\min\left\{n\ge0 \colon
\textrm{$(f^R)^n(x)$ and $(f^R)^n(y)$ lie in distinct }\Lambda_{i}
\right\}. $$
The last two properties involve information on the action of $f^R$ on unstable leaves.

{\em
\begin{enumerate}
\item[\bf(P$_3$)] {\em Regularity of   the stable foliation}: given
$\gamma,\gamma'\in\Gamma^u$, we define
$\Theta\colon\gamma'\cap\Lambda\to\gamma\cap\Lambda$  by
$\Theta(x)=\gamma^s(x)\cap \gamma$. Then
\begin{enumerate}
 \item $\Theta$ is absolutely continuous and
        $$\displaystyle \frac{d(\Theta_*\leb_{\gamma'})}{d\leb_{\gamma}}(x)=
        \prod_{i=0}^\infty\frac{\det Df^u(f^i(x))}{\det
        Df^u(f^i(\Theta^{-1}(x)))};$$
  \item letting $u(x)$ denote the density in item (a),  we have
 $$\log \frac{u(x)}{u(y)}\le C\beta^{s(x,y)},\quad\text{for
 $x,y\in\gamma'\cap\Lambda$}.$$
\end{enumerate}
\end{enumerate}
}

{\em
\begin{enumerate}
 \item[\bf(P$_4$)] {\em Bounded distortion}: 
 for $\gamma\in\Gamma^u$ and
    $x,y\in\Lambda\cap\gamma$
    $$
    \log\frac{\det D(f^R)^u(x)}{\det D(f^R)^u(y)}\le
    C\beta^{s(f^R(x),f^R(y))}.
    $$
 \end{enumerate}
}

\cre The Markov property we present here is weaker than the one in
\cite{Y1}, since includes two extra assumptions: $i)$ there are at
most finitely many $i$'s with $R_i=n$ for each $n\in\NN$; $ii)$
$R_i\ge R_0$ for some $R_0>1$ depending on the constants $C$ and
$\alpha$.  These assumptions play a role in showing the existence of a spectral gap for a transfer operator associated to the dynamics.  Here we use a more probabilistic argument,  based on
\cite{Y2}, which enables us to drop those extra assumptions. In
particular, we are able to reobtain the conclusions of~\cite{Y1} under
our weaker Markov condition.\fre

\cre\label{re.hypti} We do not assume any uniform backward contraction along
unstable leaves similar to (P4)(a) in \cite{Y1}. This would be too
restrictive for our purposes, since the application we make of our
main results does not have this property. Properties (P$_3$)(b)
and (P$_4$) are new if comparing our setup to the one in
\cite{Y1}. However, 
they can be easily  obtained from (P4) and  (P5) in
\cite{Y1}; see \cite[Lemma~1]{Y1}.\fre

\subsection{Diameter control}\label{dimcontr} Consider a sequence of {\em stopping
times} defined for the points in $\Lambda$ in the following way:
 \begin{equation}\label{def.rs}
 S_0=0,\quad S_1=R\qand S_{i+1}=S_{i}+R\circ f^{S_{i}},\quad\text{for $i\ge1$}.
 \end{equation}
We also define a nested sequence $(\cp_k)_{k\ge 0} $ of partitions
of $\Lambda$.  Let $\cp_0$ be the partition of
$\Lambda$ into the subsets $\Lambda_i$. 
Given $k\ge 1$, we say that $x$ and  $y$ belong to an element of
$\cp_k$, if both  $f^{R}(x)$ and $f^{R}(y)$ have the same stopping
times
$S_0<S_1<\cdots<S_j$ up to time $k-1$, 
and $f^{S_i}(f^R(x))$ and $f^{S_i}(f^R(y))$ belong to the same
element of $\cp_0$ for each $0\le i\le j$. By construction   we have
that $S_{j+1}(f^{R}(x))=S_{j+1}(f^{R}(y))\ge k$ and
$f^{S_{j+1}}(f^R(Q))$ is a $u$-subset.

As it will become clear in the proof of  Lemma~\ref{l.shrink}, it
will be necessary to have a control on the diameter of certain iterates of
the elements in the partitions constructed above. Take  any
$k\ge 1$ and $P\in\cp_0$.  We consider separately  the cases where $k$ is bigger than $R(P)-1$ or not.
If $k> R(P)-1$, then we define
 $$\delta_k(P)=\sup_{0\le
 \ell\le R(P)-1}\left\{\diam\left(f^{\ell}(Q \cap\gamma)\right)\colon \gamma\in\Gamma^u , \,Q\in \cp_{k-R(P)+1+\ell}, \,Q\subset P\,
\right\}.$$
On the other hand, if $k\le R(P)-1$, then we define the quantities
$$\delta^0_k(P)=\sup_{0\le \ell< R(P)-k}\left\{\diam\left(f^{\ell}(P \cap\gamma)\right)\colon \gamma\in\Gamma^u\,
\right\},$$
$$\delta_k^+(P)=\sup_{R(P)-k\le
 \ell\le R(P)-1}\left\{\diam\left(f^{\ell}(Q\cap\gamma)\right)\colon \gamma\in\Gamma^u , \,Q\in \cp_{k-R(P)+1+\ell}, \,Q\subset P\,
\right\},$$ and
$$\delta_k(P)=\sup\{\delta_k^0(P), \delta_k^+(P)\}.$$
Finally we define
 \begin{equation}\label{eq.diam}
 \delta_k=\sup_{P\in\cp_0}\delta_k(P).
 \end{equation}

Though the definition of $\delta_k$  might seem somewhat technical,
this is not so hard to calculate in practice, at least for some
examples. One we have in mind is the one that we present
at Section~\ref{se.application}, for which we show  that $\delta_k$ decays polynomially fast
with $k$; see
Section~\ref{se.diameter}.

\cre\label{re.diam} The argument in Section~\ref{se.diameter} can easily be adapted to
show that {\em $\delta_k$ decays exponentially fast  with~$k$, once we know that the
diameter of the elements $\Lambda_i$ 
decay
exponentially fast with $R_i$}. This includes all the examples studied in
\cite{Y1}, since property (P4)(a) in \cite{Y1} gives the
exponential decay for the diameters of the elements in the initial
partition with
respect to the return~time.
\fre

 \cre In the
light of \cite[Definition 2.6]{ABV} one may say that {\em $\delta_k$ decays
exponentially fast with~$k$ whenever the return time $R_i$ is a {\em
hyperbolic time} for the points in $\Lambda_i$ with respect to the
derivative restricted to the tangent direction of the leaves in
$\Gamma^u$}; see \cite[Lemma 2.7]{ABV} and recall Remark~\ref{re.diam}.\fre

\subsection{Main results}

The first result we present here asserts the existence of 
SRB measures for systems  having some hyperbolic structure,
provided the return time is integrable with respect to the
conditional of the Lebesgue measure on some local unstable leaf.

\cd
We say that an $f$-invariant probability measure $\mu$ is a {\em Sinai-Ruelle-Bowen (SRB) measure} if $f$ has no zero Lyapunov exponents $\mu$ almost everywhere, and the conditional measures on local unstable manifolds are absolutely continuous with respect to the Lebesgue measures on these manifolds.
\fd

 The
proof of the next result is quite standard and may be found in
\cite{Y1}.

\cmt\label{teor} Assume that $f$ has a hyperbolic structure
$\Lambda$ such that $\leb_\gamma(\Lambda\cap\gamma)>0$ for some
$\gamma\in\Gamma^u$. If $R$ is integrable with respect to
$\leb_\gamma$, then
 $f$ has some SRB  measure \( \mu \).
\fmt

The next result shows that the decay of correlations of the SRB measure $\mu$
given by Theorem~\ref{teor} is related to the recurrence times of
the hyperbolic structure. It has been established by Young in
\cite[Theorem 2]{Y1} a version of this result for hyperbolic
structures having exponential decay of return time. The method in
\cite{Y1} is based on the existence of a spectral gap for the
transfer operator and cannot be applied in our situation. We define
the space of H\"older continuous functions with exponent $\eta>0$
 $$H_\eta=\big\{\varphi\colon M\to\RR\;\vert\,\, \exists C>0\text{ such
 that
 }|\varphi(x)-\varphi(y)|\le C\dist(x,y)^\eta,\;\forall x,y\in
 M\,\big\}.
 $$

\cmt\label{th:DC} Assume that $f$ has a hyperbolic structure
$\Lambda$ for which (P$_1$)-(P$_4$) hold, with $\gcd\{R_i\}=1$ and
$\leb_\gamma(\Lambda\cap\gamma)>0$ for some $\gamma\in\Gamma^u$.
Given $\varphi,\psi\in H_\eta$,
\begin{enumerate}
  \item if
 $\leb_\gamma\{ R>n\}\lesssim
n^{-\alpha}$ for some $\alpha>1$, then
 $\mathcal C_{n}(\varphi,\psi;\mu) \lesssim\max\{{n^{-\alpha+1
 },\delta_n^\eta}\}$;
  \item  if
 $\leb_\gamma\{ R>n\}\lesssim
 e^{-cn^\zeta}$ for some $c>0$ and $0<\zeta\le1$, then
 there exists $c'>0$ such that
 $\mathcal C_{n}(\varphi,\psi;\mu) \lesssim \max\{e^{-c'n^\zeta},\delta_n^\eta\}$.
\end{enumerate}
\fmt

As shwon in \cite[Section~4.1]{Y1}, condition $\gcd\{R_i\}=1$ can
be replaced by the assumption that $f^n$ is ergodic with respect to
$\mu$ for every $n\ge 1$. If we omit both assumptions, then the same
conclusion holds for some power of $f$. The next result gives the
{\em Central Limit Theorem} for H\"older continuous observables which are
not a coboundary with respect to the SRB measure $\mu$.

\cmt\label{th:TLC} Under the assumptions of Theorem~\ref{th:DC}, if
 $\leb_\gamma\{ R>n\}\lesssim n^{-\alpha}$ for some $\alpha>2$, then given  \( \varphi \in
H_\eta\) for which there is no $\psi\in L^2
 (\mu)$
 with  \( \varphi = \psi\circ f - \psi \)
there exists \( \sigma>0 \) such that for every interval \( J\subset
\mathbb R \),
\[
\mu\left\{x\in M: \frac{1}{\sqrt
n}\sum_{j=0}^{n-1}\left(\varphi(f^{j}(x))-\int\varphi d\mu
\right)\in J \right\}
\stackrel{n\to\infty}{\longrightarrow}\frac{1}{\sigma \sqrt{2\pi}
}\int_{J} e^{-t^{2}/ 2\sigma^{2}}dt.
\]

\fmt

\subsection{Application}\label{se.application}
\label{ss.example} We give  a diffeomorphism where
we may apply our main results and deduce that it has an SRB measure
with polynomial decay of correlations. This is obtained by
perturbing the classical solenoid map
 in the unstable direction of one fixed
point and transforming it into an indifferent fixed point.
 Let $f\colon S^1\to S^1$ be a map of degree $d\ge
2$ with the following properties:
\begin{enumerate}
  \item[(i)] $f$ is $C^2$ on $S^1\setminus\{0\}$;
  \item[(ii)]  $f$ is $C^1$ on $S^1$ and $f'>1$ on $S^1\setminus\{0\}$;
  \item[(iii)] $f(0)=0$, $f'(0)=1$, and there is $\gamma>0$ such that
    $$-xf''(x)\approx |x|^\gamma\quad\text{for all $x\neq 0$.}$$
\end{enumerate}
Consider the solid torus $M=S^1\times D^2$, where $D^2$ is the unit
disk in $\RR^2$, and define the map $g\colon M\to M$ by
 $$g(x,y,z)=\left(f(x), \frac1{10}y+\frac12\cos x,
 \frac1{10}z+\frac12\sin x\right).$$
Let $H_\eta$ be the space of H\"older continuous functions on $M$
with exponent $\eta>0$.

\cmt\label{th:example} Let $g\colon M\to M$ be as above and take
$\varphi,\psi\in H_\eta$.
\begin{enumerate}
  \item The map $g$ admits an SRB measure $\mu$ if and only if $\gamma<1$.
  \item Assume that $\gamma<1$.  Then
  \begin{enumerate}
    \item for $\eta\ge 1-\gamma$ we have $\mathcal
    D_{n}(\varphi,\psi;\mu)\lesssim n^{1-1/\gamma}$;
    \item for $\eta < 1-\gamma$ we have $\mathcal
    D_{n}(\varphi,\psi;\mu)\lesssim n^{-\eta/\gamma}$.
  \end{enumerate}
  \item If $\gamma<1/2$, then the Central Limit Theorem holds for $\varphi\in
    H_\eta$, provided there is no $\psi\in L^2
 (\mu)$
 with  \( \varphi = \psi\circ f - \psi \).
\end{enumerate}
\fmt

It is well known that for $\gamma\ge 1$ one has
$\frac1n\sum_{j=0}^{n-1}\delta_{f^j(x)}$ converging in the weak*
topology to the Dirac measure at 0 for Lebesgue almost every $x\in
S^1$; see for example \cite{HK} and \cite{Pi}. Using the fact that
we have uniform contraction in the vertical direction, it is not hard to see that
$\frac1n\sum_{j=0}^{n-1}\delta_{g^j(x,y)}$ converges in the weak*
topology to the Dirac measure at 0 for Lebesgue almost every
$(x,y)\in S^1\times D^2$. This observation justifies the ``only if" part of the theorem above.

\section{Induced schemes}

The objects that we introduce in this section have essentially
 been all presented in \cite{BY2} and put into an abstract setting in
\cite{Y1}.

\subsection{The natural measure}\label{The return map} Fix an
arbitrary $\hat\gamma\in\Gamma^u$.  Given $\gamma\in\Gamma^u$ and
$x\in\gamma\cap\Lambda$ let~$\hat x$ be the point in
$\gamma^s(x)\cap \hat\gamma$. Defining for $x\in\gamma\cap\Lambda$
$$\hat u(x)=\prod_{i=0}^\infty\frac{\det Df^u(f^i(x))}{\det
        Df^u(f^i(\hat x))}$$
we have that $\hat u$ satisfies the bounded distortion property
(P$_3$)(b). For each $\gamma\in\Gamma^u$ let $m_\gamma$ be the
measure in $\gamma$ such that
 $$\frac{dm_\gamma}{d\leb_\gamma}=\hat u\,\um_{\gamma\cap\Lambda},$$
where $\um_{\gamma\cap\Lambda}$  is the characteristic function of
the set ${\gamma\cap\Lambda}$. These measures have been defined in
such a way that if $\gamma,\gamma'\in\Gamma^u$ and $\Theta$ is
obtained by sliding along stable leaves from $\gamma\cap\Lambda$ to
$\gamma'\cap\Lambda$, then
 \begin{equation}\label{eq.teta}
\Theta_*m_{\gamma}=m_{\gamma'}.
\end{equation}
To verify this let us show that the densities of these two measures
with respect to $\leb_{\gamma}$ coincide. Take
$x\in\gamma\cap\Lambda$ and  $x'\in\gamma'\cap\Lambda$ such that
$\Theta(x)=x'$.  By
 (P$_3$)(a) one has
   $$\frac{d\Theta_*\leb_{\gamma}}{d\leb_{\gamma'}}(x')=\frac{\hat u(x')}{\hat u(x)},$$
   which implies that
  $$\frac{d\Theta_*m_{\gamma}}{d\leb_{\gamma'}}(x')={\hat u(x)}
  \frac{d\Theta_*\leb_\gamma}{d\leb_{\gamma'}}(x')={\hat u(x')}=
 \frac{dm_{\gamma'}}{d\leb_{\gamma'}}(x').$$

\cle\label{distgama}  Assuming that
$f^R(\gamma\cap\Lambda)\subset\gamma'$ for
$\gamma,\gamma'\in\Gamma^u$, let $Jf^R(x)$ denote the Jacobian of
$f^R$ with respect to the measures $m_\gamma$ and $m_{\gamma'}$.
Then
\begin{enumerate}
    \item $Jf^R(x)=Jf^R(y)$ for every $y\in\gamma^s(x)$;
    \item there is $C_1>0$ such that for every
    $x,y\in\Lambda\cap\gamma$
    $$
    \left|\frac{Jf^R(x)}{Jf^R(y)}-1\right|\le
    C_1\beta^{s(f^R(x),f^R(y))}.
    $$
\end{enumerate}
\fle
 \dem (1) For $\leb_\gamma$ almost every $x\in\gamma\cap\Lambda$ we
 have
  \begin{equation}\label{eq.jaco}
  Jf^R(x)=\left|\det D (f^R)^u(x)\right|\cdot \frac{\hat u(f^R(x))}{\hat u(x)}.
  \end{equation}
 Denoting $\varphi(x)=\log|\det Df^u(x)|$ we may write
 $$
 \begin{array}{lllll}
  \log Jf^R(x) &=& \displaystyle\sum_{i=0}^{R-1}\varphi(f^i(x))&+&
  \displaystyle\sum_{i=0}^{\infty}\left(\varphi(f^i(f^R(x)))- \varphi(f^i(\widehat{f^R(x)})\right)\\
     & &   & - & \displaystyle\sum_{i=0}^{\infty}\left(\varphi(f^i(x))- \varphi(f^i(\hat
     x)\right)\\
     &=& \displaystyle\sum_{i=0}^{R-1}\varphi(f^i(\hat x))&+&
  \displaystyle\sum_{i=0}^{\infty}\left(\varphi(f^i(f^R(\hat x)))- \varphi(f^i(\widehat{f^R(x)})\right)
 \end{array}
  $$
Thus we have shown that $Jf^R(x)$ can be expressed just in terms of
$\hat x$ and $\widehat{f^R(x)}$, which is enough for proving the
first part of the lemma.
\smallskip

 (2) It follows from \eqref{eq.jaco} that
  $$
  \log\frac{Jf^R(x)}{Jf^R(y)}=\log\frac{\det D(f^R)^u(x)}{\det
  D(f^R)^u(y)}+
  \log\frac{\hat u(f^R(x))}{\hat u(f^R(y))}+\log\frac{\hat u(y)}{\hat u(x)}
  $$
Observing that $s(x,y)> s(f^R(x),f^R(y))$ the conclusion follows
from (P$_3$)(b) and (P$_4$).
 \cqd

\subsection{A tower extension}\label{sec.tower}

We introduce a tower extension of the dynamical system $f$
restricted to $\cup_{n\ge0}f^n(\Lambda)$; note that this space is
preserved by $f$. We define a {\em tower}
 $$\Delta=\big\{ (x,\ell)\colon \text{$x
\in \Lambda$ and $0\leq \ell < R(x)$} \big\},$$
 and a {\em tower map} $F:\Delta\to\Delta$ as
 $$
 F(x,\ell)=\left\{%
\begin{array}{ll}
    (x,\ell+1), & \hbox{if $
\ell+1<R(x)$;} \\
    (f^R(x),0), & \hbox{if $\ell+1=R(x)$.} \\
\end{array}%
\right.
 $$
The {\em $\ell^\text{th}$ level of the tower} is by definition the
set
$$\Delta_\ell =\{(x,\ell)\in \Delta \}.$$
The $0^\text{th}$- level of the tower $\Delta_0$  is naturally
identified with $\Lambda$ and we shall make no distinction between
them. Under this identification it easily follows from the
definitions that $F^R=f^R$ for each $x\in\Delta_0$. Note that the
$\ell^\text{th}$ level of the tower is a copy of the set
$\{R>\ell\}\subset \Delta_0$. Also, we easily obtain a partition
$\cp$ of $\Delta_0$ into subsets $\Delta_{0,i}$, with
$\Delta_{0,i}=\Lambda_i$ for $i\ge 1$. This partition gives rise to
partitions $\Delta_{\ell,i}$ on each tower level $\ell$, considering
 $$\Delta_{\ell,i}=\{(x,\ell)\in \Delta_\ell\,\colon\, x\in \Delta_{0,i}\}. $$
Collecting all these sets we obtain  a partition
$\mathcal{Q}=\{\Delta_{\ell,i}\}_{\ell,i}$ of $\Delta$. We introduce
a sequence of partitions $(\cq_n)_{n\ge0}$ of $\Delta$ in the
following way:
 \begin{equation}\label{qn}
    \cq_0=\cq,\qand \cq_n=\vee_{i=0}^nF^{-i}\cq \quad\text{for
$n\ge0$}.
\end{equation}
We shall denote by $Q_n(x)$ the element in $\cq_n$ containing the
point $x\in \Delta$.

We define a projection map
 \begin{equation}\label{eq.pi}
 \begin{array}{rcccl}
 \pi&\colon &  \Delta & \longrightarrow &
 \bigcup_{n\ge0}f^n(\Delta_0).
 \\
& &(x,\ell)&\longmapsto & f^\ell(x)
 \end{array}
 \end{equation}
Observe  that
 $f\circ\pi =\pi\circ F.$

 \cle\label{l.shrink} There is $C_2>0$ such that for all
$k\ge 0$ and $Q\in\cq_{2k}$  $$\diam(\pi F^k(Q ))\le
C_2\max\{\beta^k,\delta_k\}.$$ \fle

\dem Take $k\ge 0$ and $Q\in\cq_{2k}$. Given $x,y\in Q$, there is
$z\in \gamma^u(x)\cap\gamma^s(y)$. Supposing that $Q\subset
\Delta_\ell$, then $y_0=\pi F^{-\ell}(y)$  and $z_0=\pi
F^{-\ell}(z)$ are both in $\Delta_0$ and they lie on the same stable
leaf. Hence
 \begin{equation*}
 \dist(\pi F^k(y),\pi F^k(z))=\dist(\pi F^{k+\ell}(y_0),\pi
 F^{k+\ell}(z_0))=\dist(f^{k+\ell}(\pi y_0),\pi
 f^{k+\ell}(\pi z_0)).
 \end{equation*}
 Using (P$_2$) we get
\begin{equation}\label{eq.ditFs}
 \dist(\pi F^k(y),\pi F^k(z))\le C\beta^{k+\ell}.
 \end{equation}

On the other hand, we have $F^k(Q)\in\cq_k$, which implies that
$F^k(x)$ and $F^k(z)$ are both in an unstable leaf of some element
of $\cq_k$. In particular, there are $P\in\cp_0$ and $\ell<R(P)$
such that that element of $\cq_k$ is in the $\ell$-th level of the
tower over $P$. Moreover, the situations considered for defining
$\delta_k(P)$ correspond precisely to the possible cases for the
elements of $\cq_k$ over $P$. Taking into account the definition of
$\pi$, this gives
 \begin{equation*}\label{eq.ditFu}
 \dist(\pi F^k(x),\pi F^k(z))\le \delta_k(P),
 \end{equation*}
which together with \eqref{eq.ditFs} gives the desired conclusion.
\cqd

Let $m$ be the measure on $\Lambda$ whose conditional measures on
$\gamma\cap\Lambda$ with $\gamma\in\Gamma^u$ are the measures
$m_\gamma$ introduced in the previous section. This measure $m$
allows us to introduce a measure on $\Delta$ that we still denote
$m$, by letting $m\vert\Delta_\ell$ be the measure induced by the
natural identification of $\Delta_\ell$ with a subset of $\Lambda$.
We let $JF$ denote the Jacobian of $F$ with respect to this measure
$m$.

\cle\label{l.cfbeta0} There is $C_F>0$ such that for all $k\ge1$ and
all
 $x,y\in \Delta$ belonging to a same element of $\cq_{k-1}$
\begin{eqnarray*} \left| \frac{J F^k(x)}{J F^k(y)}-1\right| \leq C_F
\beta^{  s( F^k(x), F^k(y))}.
\end{eqnarray*}
\fle

\dem By Lemma~\ref{distgama}  one knows that for all $i\ge 1$ and
all $ x,y\in \Delta_{0,i}$
\begin{equation}\label{jac}
\mdlo{ \frac{J F^{ R}(x)}{J  F^{ R}(y)}-1 } \leq  C_1 \beta^{s( F^{
R}(x), F^{ R}(y))}.
\end{equation}
It follows that there is a constant $C_F>0$ such that for all
$n\ge1$ and all $x,y$ belonging to a same element of
$\vee_{j=0}^{n-1}(F^R)^{-j}\mathcal P$
\begin{equation}\label{jacfrn}
\left| \frac{J(F^R)^n(x)}{J(F^R)^n(y)}-1\right| \leq C_F
\beta^{s((F^R)^n (x),(F^R)^n(y))},
\end{equation}
In fact, if  $x$ and $y$ belong to a same element of
$\vee_{j=0}^{n-1}(F^R)^{-j}\mathcal P$, then $(F^R)^j(x)$ and
$(F^R)^j(y)$ belong to a same element of $\mathcal P$ for every
$0\le j<n$. Moreover,

\begin{equation}\label{esse}
s((F^R)^j(x),(F^R)^j(y))=s((F^R)^n(x),(F^R)^n(y))+(n-j).
\end{equation}
Then
\begin{eqnarray}
\displaystyle\log\frac{J(F^R)^n(x)}{J(F^R)^n(y)} &=&
\displaystyle\sum_{j=0}^{n-1}\log\frac{JF^R((F^R)^j(x))}{JF^R((F^R)^j(y))}\nonumber\\
&\le&\displaystyle\sum_{j=0}^{n-1}C_1\beta^{s((F^R)^n(x),(F^R)^n(y))+(n-j)-1},
\textrm{\quad by \eqref{jac} and \eqref{esse}}\nonumber\\
&\displaystyle\leq& C_F\beta^{s((F^R)^n(x),(F^R)^n(y))},
\label{logo}
\end{eqnarray}
where $C_F>0$ depends only on $C_1$ and $\beta$. This implies that
\eqref{jacfrn} holds.

From \eqref{jacfrn} we easily deduce
 that for all $k\ge1$ and all
 $x,y\in\Delta$ belonging to a same element of $\cq_{k-1}$
\begin{eqnarray} \left| \frac{JF^k(x)}{JF^k(y)}-1\right| \leq C_F
\beta^{s(F^k(x),F^k(y))}. \label{jac5}
\end{eqnarray}
To see this, we consider $JF^k(x)=J(F^R)^n(x')$ and
$JF^k(y)=J(F^R)^n(y')$, where $n$ is the number of visits of $x$ and
$y$ to $\Delta_0$ prior to time $k$, and $x', y'$ are the elements
in the bottom level $\Delta_0$ corresponding $x,y$, respectively. In
this way, we have $x',y'$ belonging to a same element of
$\vee_{j=0}^{n-1}(F^R)^{-j}\mathcal P$ and $s(x,y)=s(x',y')$. Using
\eqref{jacfrn} we obtain \eqref{jac5}. \cqd

\subsection{Quotient dynamics}\label{s.quotientdynamics}

Let $\bar\Lambda=\Lambda/\sim$, where $x\sim y$ if and only if $y\in
\gamma^s(x)$. This quotient space  gives rise to a quotient tower
$\bar\Delta$ with levels $\bar \Delta_\ell=\Delta_\ell/\sim$. A
partition of $\bar\Delta$ into $\bar\Delta_{0,i}$, that we denote by
$\bar\cp$, and a sequence $\bar\cq_n$ of partitions of $\bar\Delta$
as in \eqref{qn} are defined in a natural way.

As $f^R$ takes $\gamma^s$-leaves to $\gamma^s$-leaves and $R$ has
been defined in such a way that it does not depend on the point we
take in a same stable leaf, we may assume that we have defined the
return time $\bar R\colon \bar\Delta_0\to\NN$, the tower map $ \bar
F\colon\bar\Delta\to\bar\Delta$ and the separation time $\bar
s\colon \bar\Delta_0\times \bar\Delta_0\to \NN$ naturally induced by
the corresponding ones in $\Delta_0$ and $\Delta$. It will be
convenient to have this separation time defined in the
whole~$\bar\Delta$. This may be done by taking $ \bar s(x,y)= \bar
s(x',y')$ if $x$ and $y$ belong in a same $\bar\Delta_{l,i}$, where
$x',y'$ are the corresponding elements of $\bar\Delta_{0,i}$, and
$\bar s(x,y)=0$ otherwise.


Since \eqref{eq.teta} holds, we may introduce a measure $\bar m$ on
$\bar\Delta$ whose representative on each $\gamma\in\Gamma^u$ is
$m_\gamma$. We let $J\bar F$ denote the Jacobian of $\bar F$ with
respect to this measure $\bar m$. The first item of
Lemma~\ref{distgama} shows that the Jacobian $J\bar F$ is well
defined with respect to $\bar m$. From Lemma~\ref{l.cfbeta0} we
easily obtain:

\cle\label{l.cfbeta} For all $k\ge1$ and all
 $x,y\in\bar\Delta$ belonging to a same element of $\bar\cq_{k-1}$
\begin{eqnarray*} \left| \frac{J\bar F^k(x)}{J\bar F^k(y)}-1\right| \leq C_F
\beta^{ \bar s(\bar F^k(x),\bar F^k(y))}.
\end{eqnarray*}
\fle

It will be useful to consider  $\hat R\colon
\bar\Delta\longrightarrow \NN$ defined as
$$
\hat R(x)=\min\{ n\ge 0 \colon \textrm{ $\bar
F^n(x)\in\bar\Delta_0$}\}.
$$
 Note that
$\hat R(x)=\bar R(x)$ for all $x\in \bar\Delta_0$, and
$$
\bar m\{\hat{R}>n\}=\sum_{l>n}\bar m(\bar\Delta_l)=\sum_{l>n}\bar
m\{\bar R>l\}.
$$
We introduce the  spaces of H\"older functions in $\bar\Delta$
\begin{align*}
 \hc_{\beta}
 = &\left\{ \varphi:\bar\Delta\to\mathbb{R} \mid~ \exists C_\varphi>0\textrm{
such that }\mdlo{\varphi(x)-\varphi(y)}\leq
C_\varphi\beta^{\bar s(x,y)}
\text{ for all } x,y \in \bar\Delta \right\} \\
\hc_\beta^{+}
= & \,\big\{ \varphi \in \hc_\beta
\mid~ \exists C_\varphi>0\mbox{ such that on each
$\bar\Delta_{\ell,i} $, either
 $\varphi \equiv 0$, or} \\  &
 \phantom{ggsdggggooasdsdfasdfrt}\varphi>0\textrm{ and }
\mdlo{\frac{\varphi(x)}{\varphi(y)}-1}\leq C_\varphi \beta^{\bar
s(x,y)} \text{ for all }  x,y \in \bar\Delta_{\ell,i} \big\}.
\end{align*}
The following result gives the existence of an equilibrium measure
for the tower map and some of its properties. A proof of it is given
in  \cite[Lemma 2]{Y1} and \cite[Theorem~1]{Y2}.

\ct\label{exis medida} Assume that $\bar R$ is integrable with respect to $\bar m$. Then
\begin{enumerate}
\item $\bar F$ has a unique absolutely continuous invariant probability $\bar\nu$
equivalent~to~$\bar m$; \item $d\bar\nu/d\bar m$ belongs to
$\hc_\beta^+$ and is bounded from below by some $c>0$;
\item $(\bar F,\bar\nu)$ is exact and, hence ergodic and mixing.
\end{enumerate}
\ft


The decay of correlations for the measure $\bar\nu$ has been proved
in \cite{Y1}. This occurs at the same speed that the positive
iterates under $\bar F_*$ of measures with densities in $
\hc_\beta^+$ converge to the equilibrium $\bar\nu$. This speed  is
related to the decay of $\bar m\{\bar R>n\}$, at least for some
specific~rates.



\ct\label{Limsup} For $\varphi\in \hc_\beta^+$  let $\bar \lambda$
be the measure whose density with respect to $\bar m$~is~$\varphi$.
\begin{enumerate}
 \item If $\bar m\{\bar R>n\}\le Cn^{-\zeta}$, for some $C>0$ and $\zeta>1$,
 then there is $C'>0$ such that
 $$
 \mdlo{\bar F^{n}_*\bar\lambda-\bar\nu} \le  C'n^{-\zeta+1}.
 $$
 \item If $\bar m\{\bar R>n\}\le C e^{-cn^{\eta}}$, for some $C,c>0$ and $0<\eta
\le1$, then there are $C',c'>0$ such that
 $$
 \mdlo{\bar F^{n}_*\bar\lambda-\bar\nu} \le C' e^{-c'n^{\eta}}.$$
\end{enumerate}
Moreover, $c'$ does not depend on $\varphi$ and $C'$ depends only on
$C_\varphi$.
    \ft

A version of this theorem has been proved in  \cite[Theorem 2]{Y1}
but without establishing the dependence on the constants. This plays a crucial role in our proofs of Theorem~\ref{th:DC} and Theorem~\ref{th:TLC}. We
postpone the proof of Theorem~\ref{Limsup} to Appendix~\ref{CE}.

\section{Back to the original dynamics}


Let $\pi$ be the map from  $\Delta$ to $M$ defined in \eqref{eq.pi}.
Let also $\bar\pi$ be the projection from $\Delta$ to the quotient
space $\bar\Delta$. As observed in \cite[Sections~2 \& 4]{Y1} we
have
 $\bar\nu=\bar\pi_*\nu$ and $\mu=\pi_*\nu$. Given $\varphi,\psi\in
 H_\eta$ we define
$\tilde \psi=\psi\circ\pi$ and $\tilde \varphi=\varphi\circ\pi$.

\subsection{Decay of correlations}\label{se.dc} 

For proving Theorem~\ref{th:DC} we start by noting that for
$\varphi,\psi\in
 H_\eta$ we have
 $$
 \int (\varphi \circ f^{n}) \psi d\mu
- \int \varphi d\mu \int \psi d\mu =
 \int (\tilde\varphi \circ F^{n})
\tilde\psi d\nu - \int \tilde\varphi d\nu \int \tilde\psi d\nu,
$$
which shows that it suffices to obtain  the desired conclusions for
$\mathcal C_n(\tilde\varphi,\tilde\psi;\nu)$. This will be done in
several steps,  firstly reducing it to a problem in $\bar \Delta$
and then applying Theorem~\ref{Limsup}. 

\subsubsection*{Step 1} Fix some positive integer $k\le n/4$.
Consider a discretization $\bar\varphi_k$ of $\tilde\varphi$ defined
on $\Delta$ (or $\bar\Delta$) as
$$\bar\varphi_k\vert_A= \inf\{\tilde\varphi\circ
F^k(x)\colon {x\in A}\},\quad\text{for $A\in \mathcal Q_{2k}$.}$$
We have
\begin{equation}\label{eq.d1}
|\mathcal C_n(\tilde\varphi,\tilde\psi;\nu)-\mathcal C_{n-k}
(\bar\varphi_k,\tilde\psi;\nu)|\le C_3\delta_k^{\eta}
\end{equation}
for some $C_3$ depending only on $C_\varphi$ and $\|\psi\|_\infty$.

Actually, by Lemma~\ref{l.shrink} one knows that
$|\tilde\varphi\circ F^k-\bar\varphi_k|\le
C_\varphi(C_2\delta_k)^\eta$. To be precise, one should consider the
case $\beta^k>\delta_k$, but this would  only be relevant in the
second part of Theorem~\ref{th:DC}. However, it does not play any
special role for the conclusion.

Observing that
 $\mathcal C_n(\tilde\varphi,\tilde\psi;\nu)=\mathcal C_{n-k}
(\tilde\varphi\circ F^k,\tilde\psi;\nu)$, the left hand side of
inequality \eqref{eq.d1} is
\begin{eqnarray*}
 &\le &
 \left|\int(\tilde\varphi\circ F^k-\bar\varphi_k)\circ F^{n-k}\cdot \tilde\psi d\nu\right|+
 \left|\int(\tilde\varphi\circ F^k-\bar\varphi_k) d\nu\cdot \int\tilde\psi d\nu\right|\\
 &\le&
 2C_\varphi(C_2\delta_k)^\eta \|\psi\|_\infty
\end{eqnarray*}
We just have to take $C_3= 2C_\varphi C_2^\eta \|\psi\|_\infty$.

\subsubsection*{Step 2} Consider $\bar\psi_k$ defined similarly to
$\bar\varphi_k$ above. Let $\bar\psi_k\nu$ denote the signed measure
whose density with respect to $\nu$ is $\bar\psi_k$, and let
$\tilde\psi_k$ denote the density of $F_*^k(\bar\psi_k\nu)$ with
respect to~$\nu$. Then
 \begin{equation}\label{eq.d2}
|\mathcal C_{n-k}(\bar\varphi_k,\tilde\psi;\nu)-\mathcal
C_{n-k}(\bar\varphi_k,\tilde\psi_k;\nu)|\le C_4\delta_k^\eta,
\end{equation}
for some $C_4$ depending only on $C_\psi$ and $\|\varphi\|_\infty$.

 In fact, the left hand side of \eqref{eq.d2} is
 \begin{eqnarray*}
 &\le &
 \left|\int(\bar\varphi_k\circ F^{n-k} )(\tilde\psi-\tilde\psi_k) d\nu\right|+
 \left|\int\bar\varphi_k d\nu\int(\tilde\psi-\tilde\psi_k) d\nu\right|\\
 &\le&
 2 \|\varphi\|_\infty \cdot \left|\int(\tilde\psi-\tilde\psi_k) d\nu\right|
\end{eqnarray*}
Letting $|\cdot|$ denote the total variation of a signed measure, and noting that $$F_*^k((\tilde\psi\circ F^k)\nu)=\tilde\psi\nu,$$ we have
\begin{eqnarray*}
 \left|\int(\tilde\psi-\tilde\psi_k) d\nu\right|
 &=&
  |\tilde\psi \nu-\tilde\psi_k\nu| \\
 &=&
   |F_*^k((\tilde\psi\circ F^k)\nu)-F_*^k(\bar\psi_k\nu)| \\
   &\le& |(\tilde\psi\circ F^k-\bar\psi_k)\nu| \\
   &=& \int|\tilde\psi\circ F^k-\bar\psi_k|d\nu.
\end{eqnarray*}
By Lemma~\ref{l.shrink} one has $|\tilde\psi\circ F^k-\bar\psi_k|\le
C_\psi(C_2\delta_k)^\eta$. Take $C_4= 2C_\psi C_2^\eta
\|\varphi\|_\infty$.

\subsubsection*{Step 3} Now we show that
 \begin{equation}\label{eqdigual}
 \mathcal C_{n-k}(\bar\varphi_k,\tilde\psi_k;\nu)=\mathcal C_{n}(\bar\varphi_k,\bar\psi_k;\bar\nu)
 \end{equation}
 Indeed,
 $$
 \int( \bar\varphi_k\circ F^{n-k})\tilde\psi_kd\nu=\int\bar\varphi_k d(F^{n-k}_*(\tilde\psi_k\nu))=
 \int\bar\varphi_k d(F^{n}_*(\bar\psi_k\nu)),
 $$
and since $\bar\varphi_k$ is constant on $\gamma^s$ leaves and $F$ and $\bar F$ are semi-conjugated by $\bar\pi$,
we have
$$
 \int\bar\varphi_k d(F^{n}_*(\bar\psi_k\nu))
 =\int\bar\varphi_k d(\bar\pi_*F^{n}_*(\bar\psi_k\nu))
 =\int\bar\varphi_k d(\bar F^{n}_*(\bar\psi_k\bar\nu))
 =\int(\bar\varphi_k \circ F^{n})\bar\psi_kd\bar\nu.
 $$
 Thus we have proved that
 $$
 \int( \bar\varphi_k\circ F^{n-k})\tilde\psi_kd\nu=\int(\bar\varphi_k \circ F^{n})\bar\psi_kd\bar\nu.
 $$
 On the other hand,
 $$
 \int \bar\varphi_k d\nu \cdot\int \tilde\psi_k d\nu
 =
 \int \bar\varphi_k d\bar\nu \cdot\int d(F_*^k(\bar\psi_k \nu))=
 \int \bar\varphi_k d\bar\nu \cdot\int \bar\psi_k d\bar\nu.
 $$
These last to formulas give precisely \eqref{eqdigual}.

\subsubsection*{Step 4} 
With no loss of generality we assume that $\bar\psi_k$ is not the null
function.  Taking
 $$
 b_k=\left(\int(\bar\psi_k+2\|\bar\psi_k\|_\infty)
 d\bar\nu\right)^{-1}\qand \hat\psi_k =b_k
 (\bar\psi_k+2\|\bar\psi_k\|_\infty),
 $$
we then have 
  $$
  \int\hat{\psi_k}\bar\rho d\bar m=1, \quad\text{where $\bar\rho=\frac{d\bar\nu}{d\bar
  m}$.}
  $$
  Moreover,
  $$\frac1{3\|\bar\psi_k\|_\infty}\le b_k\le \frac1{\|\bar\psi_k\|_\infty}\qand 1\le \|\hat\psi_k\|_\infty\le 3.$$
Observe that $\hat\psi_k$ is constant on elements of $\mathcal
Q_{2k}$, since  $\bar\psi_k$ has this property.
   Let $\hat\lambda_k$ be the probability measure on
   $\bar\Delta$ whose density with respect to $\bar m$ is
$\hat{\psi_k}\bar\rho$. Then,
\begin{eqnarray}
\mdlo{\int(\bar\varphi_k\circ \bar F^n) \bar\psi_k
d\bar\nu-\int\bar\varphi_k d\bar\nu\int\bar\psi_k
d\bar\nu}&=&\frac{1}{b_k}\mdlo{\int(\bar\varphi_k\circ \bar F^n)
\hat{\psi_k}
d\bar\nu-\int\bar\varphi_k d\bar\nu\int\hat{\psi_k} d\bar\nu}\nonumber\\
&\leq&
\frac{1}{b_k}\int\mdlo{\bar\varphi_k}.\mdlo{\frac{d(\bar F^n_*\hat\lambda_k)}{d\bar m}-
\bar\rho}d\bar m .\label{edq}
\end{eqnarray}
Letting $\bar\lambda_k=\bar F_*^{2k}\hat\lambda_k$, we have
 $$\frac d{d\bar m} \bar F_*^{n}\hat\lambda_k= \frac d{d\bar m} \bar F_*^{n-2k}\bar\lambda_k,$$
which together with \eqref{edq} gives
$$
\mathcal C_{n}(\bar\varphi_k,\bar\psi_k;\bar\nu)\le
\frac{1}{b_k}\|\bar\varphi_k\|_\infty\mdlo{\bar
F^{n-2k}_*\bar\lambda_k-\bar\nu}\le
3\|\psi\|_\infty\|\varphi\|_\infty\mdlo{\bar
F^{n-2k}_*\bar\lambda_k-\bar\nu}.
$$

Let $\phi_k$ represent the density of the measure $\bar\lambda_k$
with respect to $\bar m$. The next lemma shows that
$\phi_k\in\mathcal F_\beta^+$, with the constant $C_{\phi_k}$ not
depending on $\phi_k$. This is enough for  using
Theorem~\ref{Limsup} and conclude the proof of Theorem~\ref{th:DC}.
Recall that we have taken $k\le n/4$.

\cle There is $C>0$, not depending on $\phi_k$, such that
$$
\mdlo{\phi_k(\bar x)-\phi_k(\bar y)}\leq C\beta^{\bar s(\bar
x,\bar y)}, \quad \text{for all } \bar x,\bar y \in \bar\Delta .$$
\fle
 \dem Since $\bar F^{2k}_*\bar\nu=\bar\nu$  and
 $\bar\rho=d\bar \nu/d\bar m$,  we may write
 \begin{equation}\label{eq.ro}
\bar\rho(\bar x)=\sum_{Q\in \bar\cq_{2k} }\frac{\bar
\rho\big((\bar F^{2k}\vert Q)^{-1}(\bar x)\big)}{J\bar
F^{2k}\big((\bar F^{2k}\vert Q)^{-1}(\bar x)\big)}.
\end{equation}
Recall that we have by definition
 $$\phi_k=\frac {d\bar\lambda_k}{dm}=\frac d{d\bar m} \bar F_*^{2k}\hat\lambda_k
 \qand \frac {d\hat\lambda_k}{dm}=\hat\psi_k\bar\rho.$$
 Since $\hat \psi_k$ is
constant on elements of $\cq_{2k}$, we have
\begin{equation*}
\phi_k(\bar x)=
 \sum_{Q\in \bar\cq_{2k}}c_Q
 \frac{\bar \rho\big((\bar F^{2k}\vert Q)^{-1}(\bar x)\big)}
 {J\bar F^{2k}\big((\bar F^{2k}\vert Q)^{-1}(\bar x)\big)}
,
\end{equation*}
 where $c_Q$ is constant on each $Q\in\bar\cq_{2k}$.  Hence,
 \begin{eqnarray}\label{fica}
\phi_k(\bar x)-\phi_k(\bar y)=
 \sum_{Q\in \bar\cq_{2k}}c_Q
 \left(
 \frac{\bar \rho\big((\bar F^{2k}\vert Q)^{-1}(\bar x)\big)}
 {J\bar F^{2k}\big((\bar F^{2k}\vert Q)^{-1}(\bar x)\big)}
 -\frac{\bar \rho\big((\bar F^{2k}\vert Q)^{-1}(\bar y)\big)}
 {J\bar F^{2k}\big((\bar F^{2k}\vert Q)^{-1}(\bar y)\big)}
 \right).
\end{eqnarray}
Fixing $Q\in\bar\cq_{2k}$, let $\bar x',\bar y',\in Q$ be such
$\bar F^{2k}(\bar x')=\bar x$ and  $\bar F^{2k}(\bar x')=\bar x$.
We have
 \begin{equation}\label{eq.outra}
\frac{\bar \rho(\bar x')}
 {J\bar F^{2k}(\bar x')}
 -\frac{\bar \rho(\bar y')}
 {J\bar F^{2k}(\bar y')}
 =
 \left(\frac{\bar \rho(\bar y')}
 {J\bar F^{2k}(\bar y')}\right)
 \left(
 \frac{\bar \rho(\bar x')}
 { \bar \rho(\bar y')}\frac{J\bar F^{2k}(\bar y')}
 {J\bar F^{2k}(\bar x')}
 -1
 \right).
\end{equation}
It follows from Theorem~\ref{exis medida} that there
is $C_{\bar\rho}>0$ such that 
 $$
 \log \left|
 \frac{\bar \rho(\bar x')}
 { \bar \rho(\bar y')}
 \right|\le C_{\bar\rho}\beta^{s(\bar x',\bar y')}.
 $$
On the other hand, by Lemma~\ref{l.cfbeta} there is $C_{\bar F}>0$
such that
$$ \log \left| \frac{J\bar F^{2k}(\bar y')}
 {J\bar F^{2k}(\bar x')}\right| \leq C_{\bar F}\beta^{ \bar s(\bar
F^{2k}(\bar x'),F^{2k}(\bar y'))}.
 $$
Since $s(\bar x',\bar y')\ge s(\bar F^{2k}(\bar x'),F^{2k}(\bar
y'))=s(\bar x,\bar y)$, we have
 \begin{equation}\label{label.eq}
\log
 \left|
 \frac{\bar \rho(\bar x')}
 { \bar \rho(\bar y')}\frac{J\bar F^{2k}(\bar y')}
 {J\bar F^{2k}(\bar x')}
 \right|
 \leq
 (C_{\bar\rho}+C_{\bar F})\beta^{ \bar s(\bar x,\bar y)}.
\end{equation}
Recalling \eqref{eq.ro} and the fact that $|c_Q|\le
\|\hat\psi_k\|_\infty\le 3$, it follows from \eqref{fica},
\eqref{eq.outra} and \eqref{label.eq} that there is some constant
$C>0$ not depending on $\phi_k$ such that
 \begin{eqnarray*}
\mdlo{\phi_k(\bar x)-\phi_k(\bar y)} &\le & C
 \beta^{ \bar s(\bar x,\bar y)}.
\end{eqnarray*}
Actually, we may take $C=3\|\bar\rho\|_\infty(C_{\bar\rho}+C_{\bar
F})$. \cqd


\subsection{Central Limit Theorem}

 Let $\varphi\in H_\eta$ and consider its lift to $\Delta$ defined
as $\tilde\varphi=\varphi\circ \pi$. Similarly to what we have done
at the beginning of Section~\ref{se.dc}, we easily see that for
proving Theorem~\ref{th:TLC} it is enough to obtain the Central
Limit Theorem for $\tilde\varphi$ with respect to $\nu$ on $\Delta$.
As in the study of the correlations decay, the proof uses results
from the quotient dynamics $\bar F\colon\bar\Delta\to\bar\Delta$.
Let $\bar\cb$ be the Borel $\sigma$-algebra on~$\bar\Delta$. Define
 $$\cb_0=\{\bar\pi^{-1}\bar A\colon \bar A\in\bar\cb\}\qand
 \bar\varphi_0=E_{\nu}(\tilde\varphi\mid\cb_0).
 $$
Putting together the information from \cite[Section~5.1.B]{Y1} and
Claim 1 in \cite[Section~5.2]{Y1} we easily see that it is enough to
show that
 $$
 \sum_{j\ge0}\int|P^j(\bar\varphi_0\bar\rho)|d\bar m<\infty,
 $$
where $P$ is the transfer operator associated to $(\bar F,\bar\nu)$.
The proof of the Sublemma in \cite[Section 5.2]{Y1} gives that
$\bar\varphi_0\in \hc_\beta$. Thus, if we consider $\bar\lambda$ the
measure whose density with respect to $\bar m$
is~$\bar\varphi_0\bar\rho$, then $P^j(\bar\varphi_0\bar\rho)$ is by
definition the density of $F^j_*\bar\lambda$ with respect to $\bar
m$. Hence, we just have to show that
 $$
 \sum_{j\ge0}\int\left|\frac{d}{d\bar m}\bar
F^j_*\bar\lambda\right|d\bar m<\infty.
 $$
First we ``renormalize'' $\bar\lambda$.
 Let
 $$
 b=\left(\int(\bar\varphi_0+\|\bar\varphi_0\|_\infty)
 d\bar\nu\right)^{-1}\qand \hat\varphi_0 =b
 (\bar\varphi_0+2\|\bar\varphi_0\|_\infty),
 $$
and consider $\hat\lambda$ the probability measure whose density
with respect to $\bar m $ is $\hat\varphi_0 \bar\rho$. We have
 \begin{equation}\label{int1}
 \int\hat\varphi_0\bar d\bar \nu=\int\hat\varphi_0\bar\rho d\bar
 m=1.
 \end{equation}
 Recalling that
$$\int\bar\varphi_0\bar d\bar \nu=\int\bar\varphi_0\bar\rho d\bar
m=0,$$
  we may write
 \begin{eqnarray*}
 \int\mdlo{\frac{d}{d\bar m}\bar F^j_*\bar\lambda}d\bar m
 &= &
 \int\mdlo{\frac{d}{d\bar m}\bar
F^j_*\bar\lambda- \bar\rho\int\bar\varphi_0 d\bar\nu}d\bar m.
 \end{eqnarray*}
 Using \eqref{int1} and the fact that $$\frac{d}{d\bar m}\bar
 F_*\bar\nu=\frac{d}{d\bar m}\bar\nu=\bar\rho,$$
 we obtain
 \begin{eqnarray*}
 \int\mdlo{\frac{d}{d\bar m}\bar F^j_*\bar\lambda}d\bar m &=& \frac1b
\int\mdlo{\frac{d}{d\bar m}\bar F^j_*\hat\lambda
-2\bar\rho\,\|\bar\varphi_0\|_\infty- \bar\rho\int\hat\varphi_0
d\bar\nu + 2\bar\rho\,\|\bar\varphi_0\|_\infty }d\bar m \\
 &=& \frac{1}{b}
 \int\mdlo{\frac{d}{d\bar m}\bar F^j_*\hat\lambda
- \bar\rho }d\bar m \\
 &=&
 \frac{1}{b}
 \mdlo{\bar F^j_*\hat\lambda
- \bar\nu}.
   \end{eqnarray*}
Under the hypotheses of Theorem~\ref{th:TLC} this last quantity is
clearly summable, by Theorem~\ref{Limsup}.

\section{A solenoid with intermittency}

Here we construct a  hyperbolic structure for the map $g$ defined in
Section~\ref{ss.example} which satisfies the assumptions of our main
theorems. Concerning (P$_1$)-(P$_4$), we just have to show that
(P$_1$) and (P$_4$) hold, since (P$_2$) and (P$_3$) are
trivially satisfied due to the uniform contraction of $g$  in the
vertical direction and the skew-product form of $g$. We also
need to give suitable estimates for the decay of return times and
the diameters in \eqref{eq.diam}. The conclusions of
Theorem~\ref{th:example} are then a consequence of our main results.

 The map $g$ possesses an attractor in $M$ which is precisely
  $$\Sigma=\bigcap_{n\ge0}g^n(M).$$
$\Sigma$ is locally a product of an interval by a  Cantor set. Topologically this set coincides
with the solenoid attractor for the classical case where $f$ is taken uniformly expanding in~$S^1$.

For defining the hyperbolic structure we are going to construct a
(mod 0) countable partition $\cp_0$ of an interval $I_1\subset S^1$
and associate to each element of $\cp_0$  a suitable return time
$R^*$ with respect to the map $f$. Then we take $$\Lambda=\Sigma\cap (I_1\times
D^2).$$
For each $(x,y)\in\Lambda$ we define
  $\gamma^s(x,y)=\{(x,y)\colon y\in D^2\}$ and $\gamma^u(x,y)$
as the connected component  of $\Lambda$ that contains $(x,y)$.
The $s$-subsets are precisely the sets $\Sigma\cap (P\times
D^2)$ with
$P\in\cp_0$ and the return times are taken accordingly.

\subsection{Partition and return times} Here we  recall some objects and results
from \cite[Section~6]{Y2} related to the map $f$. Let
$I_1,\dots,I_d$ be the partition of $S^1$ made by the fundamental
domains of $f$ arranged in a natural order, and assume for
definiteness that $0$ is the common endpoint of $I_1$ and $I_d$.
Letting $x_0$ be the other endpoint of $I_1$ we define a sequence
$(x_n)_n$ in $I_1$ with the property that $f(x_{n+1})=x_n$ for $n\ge
0$. Likewise, we consider $x_0'$ the endpoint of $I_d$ distinct from
$0$ and define a sequence $(x_n')_n$
in $I_d$ so that $f(x_{n+1}')=x_n'$ for $n\ge 0$. 

Let $J_n=[x_{n+1},x_n]$ and $ J_n'=[x_n',x_{n+1}']$ for $n\ge0$.
Consider the (mod 0) partition of $S^1$
 $$\ca=\{I_2,\dots,I_{d-1}\}\cup\{ J_n,J_n'; n\ge0\}.$$
Let $R=1$ on $I_2\cup\cdots I_{d-1}\cup J_0\cup J'_0$ and let
$R\vert J_n=R\vert J'_n=n+1$ for $n\ge 1$. We have $f^R(I_j)=S^1$
for $2\le j\le d-1$ and the $f^R$ images of all other elements of
$\ca$ are either $I_1\cup\cdots\cup I_{d-1}$ or $I_2\cup\cdots\cup
I_{d}$. The following results were proved in \cite[Sections 6.2 \&
6.3]{Y2}:
\begin{enumerate}
  \item Tail decay: {\em $\leb\{R>n\} \approx n^{-{1/\gamma}}$;}
  \item  Expansion: {\em there is $0<\beta<1$ such
that $ (f^R)'(x)\ge \beta^{-1}$ for every $x\in S^1\setminus\{0\}$;}
  \item Bounded distortion: {\em there is $C>0$ such that for every  $1\le i\le n$
and $x,y\in J_n$
 $$
 \log\frac{(f^i)'(x)}{(f^i)'(y)}\le
 C\frac{|f^i(x)-f^i(y)|}{|J_{n-i}|}.
 $$}
\end{enumerate}

%
%
%
This function $R$ does not qualify as a return time for  a hyperbolic structure  of $f$ satisfying
the Markov property. That role will be played by the function $R^*$ we introduce below. We use the time function $R$ to define a sequence of stopping times
$(S_i)_i$ as in~\eqref{def.rs}. We also define the sequence of
return times
 $$r_1=S_1=R,\qand r_{i+1}=S_{i+1}-S_i,\quad\text{for $i\ge1$}.$$
Using this sequence of stopping times we define the first return
time $R^*$ to $I_1$ as follows. We simply take $R^*(x)=S_i(x)$, where $i\ge 1 $
is the minimum such that $ f^{S_i}(x)\in
 I_1$.
 As shown in \cite[Section~6.2]{Y2} we have
  \begin{equation}\label{est.tail}
  \leb\{R^*>n\}\lesssim n^{-{1/\gamma}}.
  \end{equation}

Let $\cp_0$ be the Markov partition  of $I_1$ associated to $R^*$. Naturally associating the return times to the $s$-subsets described above, then \eqref{est.tail} gives the tail estimate that we need. Property
(P$_4$) is an easy consequence of the bounded distortion and
expansion above, since the estimates on the derivative of $g$ in the unstable direction are given by~$f$.

\subsection{Diameter estimate}\label{se.diameter}  Now we are going to show that  $\delta_k\lesssim
{1}/{k^{1/\gamma}}$, where
$\delta_k$ is the quantity defined in \eqref{eq.diam}. Taking into account  the uniform contraction on the stable direction, we just have to obtain the desired control on the unstable one.  We start by proving the following auxiliary result.

\cle\label{le.sN} Let $X$ be an interval in $S^1$ whose points have
the same stopping times $S_1,\dots,S_N$, for some $N\ge 1$, with
$S_N\ge k$, and such that $f^{S_N}(X)=I_1$. Then
$|X|\lesssim1/k^{1/\gamma}. $
 \fle

 \dem
Let $r_1,\dots,r_{N}$ be the return times of points in $X$.  Since
we are assuming that  $S_N=r_1+\cdots+r_N\ge k$, there must be some
$1\le m\le N$ such that $r_m\ge k/N$. Let
$$Y=f^{r_1+\cdots +r_{m-1}}(X).$$ 
Considering the interval $I\in \ca$ such that $Y\subset I$ we have
$R\vert I=r_m$. 
Bounded distortion yields
 $$|Y|\lesssim \frac{|f^{r_m}(Y)|}{|f^{r_m}(I)|}\cdot |I|.$$
On the other hand, the tail decay and expansion estimates
 give
  $$
  |I|\lesssim\left( \frac1{{r_m}}\right)^{1/\gamma},\quad |X|\le \beta^{m-1}|Y| \qand |f^{r_m}(Y)|\le \beta^{N-m}|I_1|.
  $$
 Taking into account the choice of $m$ we obtain
 $$|X|\lesssim \beta^N\left(\frac{N}{k}\right)^{1/\gamma} \lesssim
 \frac{1}{k^{1/\gamma}},$$
 and so we are done.
 \cqd

Take
$P\in\cp_0$ and $k\ge1$. We consider the three possible cases of
sets whose diameters have to be controlled. The first two correspond to $k\le R^*(P)-1$, and the last one corresponds to $k> R^*(P)-1$; recall the definition of $\delta_k$ in Section~\ref{dimcontr}.

\subsubsection*{Case 1}

Assume first that $k\le R^*(P)-1$ and $0\le\ell<R^*(P)-k$. There is
$m\ge 1$ such that $R^*(P)=r_1+\cdots +r_m$. Considering $r_0=0$,
let $0\le p < m$ be such that
 $$ r_0+\cdots +r_p\le \ell < r_0+\cdots +r_{p+1}.$$
 Letting $\ell'=\ell-(r_0+\cdots +r_p)$, we have
 $$
 \ell'<r_{p+1}\qand \ell'<r_{p+1}+\cdots +r_m-k.
 $$
Now, if $i+1=m$, then
  $$|f^\ell(P)|=|f^{\ell'}(f^{r_0+\cdots +r_p}(P))|\lesssim \left(
  \frac1{{r_m-\ell'}}\right)^{1/\gamma}\le \frac{1}{k^{1/\gamma}}.$$
Otherwise, for $p+1<m$ we use Lemma~\ref{le.sN} with
$X=f^{r_0+\cdots +r_p+\ell'}(P)$, $N=m-p$,
  and
  $$S_j=r_{p+1}+\cdots +r_{p+j}-\ell', \quad\text{for $1\le j\le N$},$$
thus obtaining
$$|f^\ell(P)|=|X|\lesssim  \frac{1}{k^{1/\gamma}}.$$ Observe that
$r_{p+1}-\ell'$ is still a return time, which then implies that
$S_1$ is well defined.

\subsubsection*{Case 2}  Assume now that $k\le R^*(P)-1$ and $R^*(P)-k\le
\ell\le R^*(P)-1$. We simply write $R^*$ for $R^*(P)$. Take $Q\in
\cp_{k-R^*+1+\ell}$ with $Q\subset P$. By construction, there is $j
\ge 0$ such that points in $f^{R^*}(Q)$ have the same stopping times
$S_{1}^*,\dots,S_{j}^*$ up to time $k-R^*+\ell$. Moreover,
$f^{S_{j+1}^*}(f^{R^*}(Q))=I_1$ and
\begin{equation}\label{eq.sjdes}
 S_{j+1}^*\ge k-R^*+1+\ell.
\end{equation}
There are integers $m,n\ge1$ and return times $r_1,\dots,r_{m+n}$
such that
\begin{equation}\label{eq.rest}
 R^*=r_1+\cdots +r_m \qand S_{j+1}^*=r_{m+1}+\cdots
 +r_{m+n}.
\end{equation}
 It follows from \eqref{eq.sjdes} and \eqref{eq.rest} that
 $$
 \ell <r_{1}+\cdots
 +r_{m+n}-k.
 $$
Considering $0\le p < m$ such that
 $$ r_0+\cdots +r_p\le \ell < r_0+\cdots +r_{p+1},$$
 where $r_0=0$ as before, and taking $\ell'=\ell-(r_0+\cdots +r_p)$, we have
 $$
 \ell'<r_{p+1}\qand \ell'<r_{p+1}+\cdots +r_{m+n}-k.
 $$
The proof now follows as in the previous case.

\subsubsection*{Case 3}  The case $k> R(P)-1\ge\ell$
is treated as Case 2.

\appendix\section{Mixing rates for tower maps}\label{CE}

The goal of this section is to prove Theorem~\ref{Limsup}. We follow the scheme of \cite{Y2} with a delicate control on the constants. The only exception is Subsection~\ref{se.gouezel} where we use results from~\cite{G}. The
setting will be the same of Subsection~\ref{s.quotientdynamics}. For
the sake of notational simplicity we shall drop all bars. 

 Let $\lambda$
and $\lambda'$ be probability measures in  $\Delta$ whose densities with respect to $m$
belong to~$\hc_\beta^+$. Let
 $$\varphi=\frac{d\lambda}{dm}\qand \varphi'=\frac{d\lambda'}{dm},$$
 and consider $C_\varphi$, $C_{\varphi'}$ as in the definition of $\hc_\beta^+$.

\subsection{Main estimates}\label{RS}

Consider the product map $F\times F:\Delta \times \Delta \to\Delta
\times \Delta$, and $P=\lambda \times \lambda'$ the product measure
on $\Delta\times\Delta$. Let $\pi$, $\pi':\Delta \times\Delta \to
\Delta$ be the projections on the first and second coordinates
respectively. Note that $F^n \circ \pi = \pi \circ (F \times F)^n$.
Consider the partition  $\cq:=\{\Delta_{l,i}\}$ of $\Delta$, and the
partition  $\cq \times \cq$ of $\Delta \times \Delta$. Note that
each element of $\cq \times \cq$ is sent bijectively by $F\times F$
onto a union of elements of $\cq \times \cq$. For each $n\ge 1$, let
$$
(\cq\times\cq)_n:= \bigvee_{i=0}^{n-1}(F\times
F)^{-i}(\cq\times\cq),
$$
and let $(\cq\times\cq)_n(x,x')$ be the element of
$(\cq\times\cq)_n$ that contains $(x,x')\in \Delta\times\Delta$.

Since $(F,\nu)$ is mixing and the density of $\nu$ with respect to  $m$ belongs to $L^\infty(m)$,
we may find $n_0\in \mathbb{N}$ and $\gamma_0>0$ such that $m(F^{-n}(\Delta_0)\cap\Delta_0)\geq \gamma_0$
for all $n\geq n_0$.
 Then we introduce a sequence of
{\em stopping times} $0\equiv \tau_0<\tau_1<\tau_2<...$  in
$\Delta\times\Delta$ given by
\begin{eqnarray*}
\tau_1(x,x')&=&n_0+\hat{R}(F^{n_0}(x)),\\
\tau_2(x,x')&=&\tau_1 + n_0+\hat{R}(F^{\tau_1}(x')),\\
\tau_3(x,x')&=&\tau_2 + n_0+\hat{R}(F^{\tau_2}(x)),\\
\tau_4(x,x')&=&\tau_3 + n_0+\hat{R}(F^{\tau_3}(x')),\\
&\vdots&
\end{eqnarray*}
with the falls to the ground level $\Delta_0$ alternating between $x$ e $x'$. This implies that
$\tau_{i+1}-\tau_{i}\ge n_0$  for all $i\ge 1$. 
We define the {\em simultaneous return time} $T:\Delta\times\Delta\to\mathbb{N}$ as
$$
T(x,x')=\min\left\{\tau_i \colon  (F^{\tau_i}(x),F^{\tau_i}(x'))\in\Delta_0\times\Delta_0, \text { with $i\geq 2$}\right\}.
$$
Note that we have  $T\ge 2n_0$. Since $(F,\nu)$ is mixing, then
$(F\times F,\nu\times\nu)$ is ergodic, and so $T$ is well-defined $m\times m$ almost everywhere. Observe that if $T(x,x')=n$,
then
$$
T\vert_{(\cq\times\cq)_n(x,x')}\equiv n\quad \text{and}
\quad(F\times F)^n((\cq\times\cq)_n(x,x'))=\Delta_0\times \Delta_0.
$$

Now we define a  sequence  $\xi_1<\xi_2<\xi_3<...$ of partitions of
$\Delta\times \Delta$. First we take $
\xi_1(x,x')=(F^{-\tau_1(x)+1}\cq)(x)\times\Delta. $ The partition
$\xi_1$ is formed by sets of the form $\Gamma=A\times\Delta$ where
$\tau_1$ is constant on $\Gamma$ and $F^{\tau_1}$ sends $A$
bijectively to $\Delta_0$. For $i>1$, if $i$ if even (resp. odd), we
define $\xi_i$ as the refinement of $\xi_{i-1}$ obtained by
partitioning $\Gamma\in\xi_{i-1}$ in the
 $x'$ direction (resp. $x$ direction) into sets $\tilde{\Gamma}$ such that $\tau_i$ is constant on each
$\tilde{\Gamma}$ and $F^{\tau_i}$ sends $\pi'(\tilde{\Gamma})$ (resp. $\pi(\tilde{\Gamma})$) bijectively
to $\Delta_0$. 
It will be useful to consider $\xi_0=\{\Delta\times\Delta\}.$ Let us
mention two useful properties about the measurability of the
functions  with respect to the partitions defined above:
\begin{itemize}
{\em
    \item $\tau_1, \tau_2, ..., \tau_i$ are $\xi_i$-measurable for each $i\ge 1$;
    \item $\{T=\tau_i\}$ and $\{T>\tau_i\}$ are  $\xi_{i+1}$-measurable for each $i\ge 1$.
}
\end{itemize}
This follows from the construction of the objects. Now we present
the main estimates we need on $\{\tau_i\}$ and $T$, whose proofs we
postpone to Section~\ref{simultaneous}.
\begin{enumerate}
\item[(E$_1$)] \label{A} {\em There is $\varepsilon_0=\varepsilon_0(C_\varphi, C_{\varphi'})>0$ such that $
 P\{T=\tau_i \mid \Gamma\} \geq \varepsilon_0
 $ for
$i\geq 2$ and $\Gamma \in \xi_i$ with $T \mid { \Gamma }
> \tau_{i-1}$.
 The dependence of $\vare_0$ on $C_\varphi$ and $C_{\varphi'}$ can
 be removed if we consider $i\ge i_0(C_\varphi, C_{\varphi'})$.
 }

\smallskip

\item[(E$_2$)]\label{B}{\em There is $K_0=K_0(C_\varphi, C_{\varphi'})>0$ such that $
P\{\tau_{i+1}-\tau_i>n_0 + n \mid\Gamma\} \leq K_0m\{\hat{R}>n\}
$ for  $i\geq 0$,  $\Gamma
\in \xi_i$ and $n\geq 0$.
The dependence of $K_0$ on $C_\varphi$ and $C_{\varphi'}$ can
 be removed if we consider $i\ge i_0(C_\varphi, C_{\varphi'})$.
}
\end{enumerate}
Let  $0\equiv T_0<T_1<T_2<\cdots $ be {\em stopping times}  in
$\Delta\times\Delta$ given by
\begin{equation}\label{tt}
    T_1=T ,\quad\text{ and }\quad T_n=T_{n-1} + T\circ (F\times F)^{T_{n-1}},\quad \text{for } n\ge2.
\end{equation}

\begin{enumerate}

\item[(E$_3$)] \label{C}{\em There are $K_1=K_1(C_\varphi, C_{\varphi'})>0$ and $\varepsilon_1>0$
(not depending  on $\varphi$ or $\varphi'$) such that $\mdlo{F_*^n\lambda-F_*^n\lambda'}\leq 2P\{T>n\}+K_1 \sum_{i=1}^{\infty}(1-\varepsilon_1)^i P\{T_i\leq n <
T_{i+1}\}$ for $n\ge 1$.
 }
\smallskip

\item[(E$_4$)] \label{D}{\em There is $K_2=K_2(C_\varphi, C_{\varphi'})>0$ such that $P\{T_{i+1}-T_i>n 
\} \leq K_2 (m\times m)\{T>n\}$ for  $i\geq 0$. 

}
\end{enumerate}

\subsection{Convergence to the equilibrium}
We shall use (E$_1$)-(E$_4$) to prove Theorem~\ref{Limsup}.
Let~$\nu$ be the measure given by Theorem~\ref{exis medida}. Observe
that $\nu$ is a fixed point for~$F_*$, whose density with respect to
$m$ belongs to $\mathcal F^+_\beta$.  Theorem~\ref{Limsup} follows
just by taking $\lambda'=\nu$, once we obtain the upper
bound  for  $\mdlo{F_*^n\lambda-F_*^n\lambda'}$. 

We start by observing that for each $i\ge 1$ we have
\begin{equation}\label{eq.start}
P\{T_{i}\leq n < T_{i+1}\}\leq
\sum_{j=0}^iP\ch{T_{j+1}-T_{j}>\frac{n}{i+1}}.
\end{equation}
Actually, since we have $T_{i+1}>n$, there must be some  $0\leq
j\leq i$ with $T_{j+1}-T_{j}> n / (i+1) $. For otherwise
$$
T_{i+1}=\sum_{j=0}^{i}(T_{j+1}-T_{j})\leq\sum_{j=0}^{i}\frac{n}{i+1}=n,
$$
which is an absurd. Hence
 $$
 \{T_{i}\leq n < T_{i+1}\}\subset \bigcup_{j=0}^i\ch{T_{j+1}-T_{j}>\frac{n}{i+1}},
$$
which gives \eqref{eq.start}. It follows respectively from (E$_3$),
\eqref{eq.start} and (E$_4$) that
\begin{equation*}
\begin{array} {llll}
\mdlo{F_*^n\lambda-F_*^n\lambda'}
 & \leq&2P\{T>n\}+
\displaystyle K_1\sum_{i=1}^\infty (1-\varepsilon_1)^i P\ch{T_{i}\leq n < T_{i+1}}, & \\
& \leq&2P\{T>n\}+\displaystyle
  K_1\sum_{i=1}^\infty
(1-\varepsilon_1)^i\sum_{j=0}^iP\ch{T_{j+1}-T_{j}>\frac{n}{i+1}},&\\
&\leq&2P\{T>n\}+\displaystyle
K_1K_2\sum_{i=1}^\infty(1-\varepsilon_1)^i(i+1)(m\times
m)\ch{T>\frac{n}{i+1}}.&
\end{array}
\end{equation*}
Observe that both in the polynomial and stretched exponential cases,
as long we obtain the desired decay for $P\{T>n\}$, then taking
$P=m\times m$ it immediately follows that
$$\sum_{i=1}^\infty(1-\varepsilon_1)^i(i+1)(m\times
m)\ch{T>\frac{n}{i+1}}$$ decays at the same speed of $P\{T>n\}$.
Consequently, we are left to  estimate  $P\{T>n\}$. At this point we
distinguish the polynomial and stretched exponential cases.

\subsubsection{Polynomial decay} Assume there are $C>0$ and $\alpha>1$ such that $m\{R>n\}\leq {C}{n^{-\alpha}}$ for all
$n\ge1$. Then, there is $\hat C>0$ (depending only on $C$ and
$\alpha$) such that
\begin{equation}\label{e.hatr}
m\{\hat{R}>n\}=\sum_{l>n}m\{R>l\}\leq \hat C{n^{-\alpha+1}}.
\end{equation}
Recall that $T\ge 2n_0$ by construction. We write
\begin{eqnarray}\label{sn}
P\{T>n\}=\sum_{1\leq i < \left[\frac{n}{2n_0}\right]}P\{T>n\colon
\tau_{i}\leq n < \tau_{i+1}\}+ P
\left\{T>n\colon\tau_{\left[\frac{n}{2n_0}\right]}\leq n\right\}.
\end{eqnarray}
Since $\{T>\tau_{i-1}\}$ is $\xi_i$-measurable, conditioning on the
elements of the partition $\xi_i$ and using  (E$_1$) it yields for
$i\ge 2$
 \begin{equation}\label{eqest1}
 P\ch{T>\tau_i\mid {T > \tau_{i-1}}}=1-P\ch{T=\tau_i \mid {T > \tau_{i-1}}}\ge
 1-\vare_0.
 \end{equation}
From \eqref{eqest1} we obtain for  $n\ge 4n_0$
\begin{eqnarray}
P \left\{T>n\colon \tau_{\left[\frac{n}{2n_0}\right]}\leq
n\right\}&\leq&
P\ch{T>\tau_{\left[\frac{n}{2n_0}\right]}}\nonumber\\
&=&P\ch{T>\tau_1}\cdot
\prod_{i=2}^{\left[\frac{n}{2n_0}\right]}P\ch{T>\tau_i
\mid {T > \tau_{i-1}}}\nonumber\\
&\leq& (1-\varepsilon_0)^{\left[\frac{n}{2n_0}\right]-1}.
\label{seco}
\end{eqnarray}
Since the dependence of $\vare_0$ on $P$ can be removed if we
consider $i\ge i_0$ for some $i_0=i_0(P)$,  we are left to
compute the decay of 
\begin{equation}\label{eq.sum2}
   \sum_{1\leq i < \left[\frac{n}{2n_0}\right]}P\{T>n\colon
\tau_{i}\leq n < \tau_{i+1}\}.
\end{equation}
For each $i\ge 1$ we have
$
P\ch{T>n \colon \tau_{i}\leq n < \tau_{i+1}}   \leq  P\ch{T>\tau_{i}
\colon n<\tau_{i+1}}
$.
 As in \eqref{eq.start}
we may show that
 $$
 \ch{T>\tau_{i} \colon n<\tau_{i+1}}\subset \bigcup_{j=0}^i\ch{T>\tau_{i} \colon \tau_{j+1} -
\tau_{j}>\frac{n}{i+1}},
$$
which then gives
\begin{equation}\label{eq.pp}
P\ch{T>n \colon \tau_{i}\leq n < \tau_{i+1}}\le \sum_{j=0}^i
P\ch{T>\tau_{i} \colon \tau_{j+1} - \tau_{j}>\frac{n}{i+1}}.
\end{equation}

Our next goal  is to estimate  the terms in the sum~\eqref{eq.pp}.
Consider first the terms with  $i,j\geq 2$. We write
\begin{equation}\label{fff}
P\ch{T>\tau_{i} \colon \tau_{j+1} - \tau_{j}>\frac{n}{i+1}}=A\cdot B\cdot C,
\end{equation}
with
\begin{align*}
A&=P\{T>\tau_1\} \cdot \prod_{k=2}^{j-1} P\{T>\tau_k\mid T>\tau_{k-1}\}\\
B&=P\ch{T>\tau_{j}\colon \tau_{j+1} -
\tau_{j}>\frac{n}{i+1}\mid T>\tau_{j-1}};\\
C&= \prod_{k=j+1}^i P\ch{T>\tau_{k}\mid T>\tau_{k-1}\colon \tau_{j+1} - \tau_{j}>\frac{n}{i+1}}. 
\end{align*}
Observe that $A=P\{T>\tau_1\}$ when $j=2$, and $C$ is void when $j=i$.
 Arguing as in \eqref{eqest1}, from estimate (E$_1$) one gets
\begin{equation}\label{eqA}
    A\leq (1-\ve_0)^{j-2}.
\end{equation}
Conditioning on $\xi_k$ and using (E$_1$), we have that
each term in C is also bounded
from above by
$1-\ve_0$, which then gives
\begin{equation}\label{eqC}
    C\leq(1-\ve_0)^{i-j}.
\end{equation}
Since $\{T>\tau_{i-1}\}$ is $\xi_i$-measurable, conditioning on
elements of $\xi_i$ and using  (E$_2$) we get
\begin{equation}\label{eqB}
    B\leq P\ch{\tau_{j+1}-\tau_j>\frac{n}{i+1}\mid T>\tau_{j-1}}\le K_0m\ch{\hat{R}>\frac{n}{i+1}-n_0}.
\end{equation}
Using \eqref{e.hatr} and the fact that $i < \left[\frac{n}{2n_0}\right]$ we obtain
\begin{equation*}\label{tti}
B  \leq K_0\hat C \left(\frac{n}{i+1}-n_0\right )^{-\alpha+1} \le
K_0\hat C 2^{1-\alpha}\left(\frac{n}{i+1}\right )^{-\alpha+1}
\end{equation*}
From \eqref{fff}, \eqref{eqA}, \eqref{eqC} and \eqref{eqB} we deduce for $i,j\ge2$ 
\begin{equation}\label{eqlarge}
P\ch{T>\tau_{i} \colon \tau_{j+1} - \tau_{j}>\frac{n}{i+1}}
 \leq
K_0\hat C 2^{1-\alpha}\left(\frac{n}{i+1}\right
)^{-\alpha+1}(1-\ve_0)^{i-2}.
\end{equation}

 Let us consider now the small terms in the sum~\eqref{eq.pp}.
  For $i\geq 2$ and $j=0,1$ we write
\begin{align*}
\lefteqn{\hspace{-8cm}P\ch{T>\tau_i\colon\tau_{j+1}-\tau_j>\frac{n}{i+1}}}\\
\lefteqn{\hspace{-6.5cm}\leq P\lefteqn{\ch{T>\tau_1\colon\tau_{j+1}-\tau_j>\frac{n}{i+1}}\cdot
\prod_{k=2}^i P\ch{T>\tau_{k}\mid T>\tau_{k-1}\colon \tau_{j+1} - \tau_{j}>\frac{n}{i+1}}
}}\\
\lefteqn{\hspace{-6.5cm}\leq
P\lefteqn{\ch{\tau_{j+1}-\tau_j>\frac{n}{i+1}}\cdot \prod_{k=2}^i
P\ch{T>\tau_{k}\mid T>\tau_{k-1}\colon \tau_{j+1} -
\tau_{j}>\frac{n}{i+1}} }}.
\end{align*}
We treat this case arguing as before, thus obtaining
\begin{equation}\label{eqsmall1}
    P\left\{T>\tau_i\colon\tau_{j+1}-\tau_j>\frac{n}{i+1}\right\}\le
    K_0\hat C 2^{1-\alpha}\left(\frac{n}{i+1}\right )^{-\alpha+1}(1-\vare_0)^{i-1}.
\end{equation}
 Finally, for $i=1$ and $j=0,1$, we have
 \begin{equation}\label{eqsmall2}
P\ch{T>\tau_{1} \colon \tau_{j+1} - \tau_{j}>\frac{n}{i+1}}
 \leq P\ch{\tau_{j+1} - \tau_{j}>\frac{n}{i+1}}\le
K_0\hat C 2^{1-\alpha}\left(\frac{n}{i+1}\right )^{-\alpha+1}. 
\end{equation}
Using \eqref{eq.pp}, \eqref{eqlarge}, \eqref{eqsmall1} and \eqref{eqsmall2}
we get
 $
P\ch{T>n : \tau_{i}\leq n < \tau_{i+1}} \leq
C_0(1-\varepsilon_0)^{i}(i+1)^{\alpha}n^{-\alpha+1},
 $
where $C_0$ is a constant depending only on $K_0$, $\hat C$,
$\alpha$ and $\vare_0$. This yields the desired bound for
\eqref{eq.sum2} in the polynomial case.


\subsubsection{Stretched exponential decay}\label{se.gouezel} Assume that there are $C,c>0$ and $0<\eta\le1$
 such that $\leb\{ R>n\}\le
 Ce^{-cn^\eta}$ for all $n\ge1$. Then there is $\hat C>0$  such that
\begin{equation*}\label{e.hatr2}
m\{\hat{R}>n\}=\sum_{l>n}m\{R>l\}\leq \hat Ce^{-cn^\eta}.
\end{equation*} The conclusion in this case is
 a consequence of (E$_1$)-(E$_2$) and the next lemma, which can
 easily
 be obtained from
 \cite[Lemma 4.2]{G} by taking $L=1$, $\tau=T$, 
 $\mu=P$ and $t_j=\tau_j$.

 \cle
Assume that there are $\varepsilon_0>0$ and $K_0>0$ such that for
all $i\geq 2$ and $\Gamma \in \xi_i$ with $T \mid { \Gamma }
> \tau_{i-1}$ we have
\begin{enumerate}
               \item $
 P\{T=\tau_i \mid \Gamma\} \geq \varepsilon_0;
 $
               \item $
P\{\tau_{i+1}-\tau_i> n \mid\Gamma\} \leq K_0e^{-cn^\eta}. $
             \end{enumerate}
Then there exist $C',c'>0$ such that $P\{T>n\}\le C'e^{-c'n^\eta}$.
 \fle

 For the sake of completeness one must verify that the constants
 $C'$ and $c'$ obey the final requirement  of Theorem~\ref{Limsup}.
 Actually, it is proved in  \cite[Lemma 4.2]{G}  that there are
a measurable function $k$   and a measurable set $B_n$ such that for
$q(n)=[\alpha n^\eta]$, with small~$\alpha>0$, we have
 $\{T>n\}\subset\{k>q(n)\}\cup B_n$
  (recall estimate (27) in \cite{G}) with
 $P\{k>q(n)\}\le (1-\vare_0)^{q(n)},$
and for some  positive integer $K$ only depending on $K_0$ and
$\eta$,
 $$
 P(B_n)\le 2^{q(n)}\sum_{p\ge n/2} C e^{-cp^\eta},\quad\text{for  $n\ge
 K$.}
 $$
Hence, taking $\alpha>0$ sufficiently small we obtain the desired
conclusion.


\subsection{Main estimates}\label{simultaneous}

Here we obtain estimates (E$_1$)-(E$_4$). We start with some
preliminary results on distortion control that will enable us to
prove (E$_1$) and (E$_2$).

\cle \label{A2} There is $C_0=C_0(C_\varphi)>0$ such that for every
$k\ge 1$  and  $A\in \vee_{i=0}^{k-1}F^{-i}\cq$ with
$F^k(A)=\Delta_0$, we have for $\mu=F_*^k(\lambda \vert A)$ and all
$x,y \in \Delta_0$
$$\left|{\frac{d\mu}{dm}(x)}\bigg/{\frac{d\mu}{dm}(y)}\right|\leq C_0.$$
Moreover, the dependence of $C_0$ on $C_\varphi$ can be removed if
we assume that the number of visits $j\le k$ of $A$ to $\Delta_0$ is
bigger than some $j_0=j_0(C_\varphi)$.
 \fle \dem Let $x_0, y_0\in A$ be such that
$F^k(x_0)=x$ and $F^k(y_0)=y$. Using \eqref{jac5} and the fact that $\varphi\in \mathcal F_\beta$ we have
\begin{equation*}
\left|{\frac{d\mu}{dm}(x)}\bigg/{\frac{d\mu}{dm}(y)}\right|
=\left|\frac{\vf(x_0)}{JF^k(x_0)}\cdot\frac{JF^k(y_0)}{\vf(y_0)}\right|
 \le   \frac{\vf(x_0)}{\vf(y_0)}\cdot \left|\frac{JF^k(y_0)}{JF^k(x_0)}\right|
 \leq (1+C_\varphi\beta^j)(1+C_F),
\end{equation*}
where $j$ is the number of visits of $A$ to $\Delta_0$ prior to $k$.
\cqd

The next  result is proved in \cite[Sublemma 1]{Y2}.

\cle \label{A1} There is $M_0>0$ such that
$\displaystyle
\frac{dF_*^nm}{dm}\leq M_0
$
for all $n\ge 1$.
\fle


\subsubsection{Proof of (E$_1$)}

Assume without loss of generality  that $i$ is even, and take  $\Gamma\in\xi_i$ as in the statement of (E$_1$).
 We have $\Gamma=A\times B$ with $A,B\subset \Delta$, where
 $A$ is sent bijectively by $F^{\tau_{i-1}}$ to $\Delta_0$ and
 $F^{\tau_{i-1}}(B)$ is contained in some $\Delta_{l,j}$. At time $\tau_i$ we have
$F^{\tau_i}(B)=\Delta_0$ and $F^{\tau_i}(A)$ is spread over several parts of
 $\cup \{\Delta_l: l\leq \tau_i - \tau_{i-1}\}$.
The set $\{T=\tau_i\}\cap \Gamma$ has the form $A'\times B$ where $A'$ is the
set of points in $A$ which are sent to $\Delta_0$ by  $F^{\tau_{i-1}}$ and return to $\Delta_0$ by
$F^{\tau_i-\tau_{i-1}}$.
Letting $\mu=F^{\tau_{i-1}}_*(\lambda\vert A)$ 
%
we may write
\begin{equation*}\label{e0}
P\{T=\tau_i\mid\Gamma\}=\frac{\lambda(A')}{\lambda(A)}=
\frac{\mu\left( F^{-(\tau_i-\tau_{i-1})}(\Delta_0)\cap\Delta_0\right)}{\mu(\Delta_0)}.
\end{equation*}
Note that Lemma~\ref{A2} applies to $\mu$, thus giving
\begin{equation*}\label{e00}
P\{T=\tau_i\mid\Gamma\}\ge C_0^{-2}
\frac{m\left( F^{-(\tau_i-\tau_{i-1})}(\Delta_0)\cap\Delta_0\right)}{m(\Delta_0)}.
\end{equation*}
Recall that  $n_0$ has been chosen in such a way that there is $\gamma_0$ such that $
m(F^{-n}(\Delta_0)\cap\Delta_0)\geq\gamma_0>0,
$ for all $n\geq n_0$.
By construction we have $\tau_i-\tau_{i-1}\geq n_0$. This is enough for concluding that
 there is some $\varepsilon_0=\varepsilon_0(C_\varphi)>0$ for which $P\{T=\tau_i\mid\Gamma\}\ge \varepsilon_0$.
The other case ($i$ odd) gives the dependence of $\varepsilon_0$
also on $C_{\varphi'}$. These dependencies can be removed if we take
$i$ large enough, according to Lemma~\ref{A2}.

\subsubsection{Proof of (E$_2$)}

For $i=0$ we have
 $$
P \{
\tau_1>n_0+n\}=(F^{n_0}_*\lambda)\{\hat{R}>n\}\leq\norm{\frac{d\lambda}{dm}}_\infty
M_0 m\{\hat{R}>n\},
 $$
and   for $i=1$
 $$
 P \{
\tau_2-\tau_1>n_0+n\}=(F^{\tau_1+n_0}_*\lambda')\{\hat{R}>n\}\leq\norm{\frac{d\lambda'}{dm}}_\infty
M_0 m\{\hat{R}>n\},$$
 which obviously give upper bounds depending on $C_\varphi$ and
 $C_{\varphi'}$.

Let us consider now the case $i\geq 2$. Assume for definiteness that
$i$ is even. Considering the probability measure
$$\mu=\frac{1}{P(\Gamma)}F^{\tau_{i-1}}_*\pi_*(P\vert\Gamma)$$ we have
\begin{eqnarray*}
P\{\tau_{i+1}-\tau_i>n_0+n\mid\Gamma\}&=&\left(
F^{(\tau_i-\tau_{i-1})+n_0}_*\mu \right)\{\hat{R}>n\}\\
&\leq&\norm{\frac{d}{dm}\left(F^{(\tau_i-\tau_{i-1})+n_0}_*\mu\right)}_\infty
m\{\hat{R}>n\}\\
&\leq& M_0\norm{\frac{d\mu}{dm}}_\infty m\{\hat{R}>n\}, \qquad \textrm{ by Lemma \ref{A1}}.
\end{eqnarray*}
Using  Lemma~\ref{A2} one has that $\norm{{d\mu}/{dm}}_\infty$ is
bounded from above by some constant only depending on $C_0$.
Moreover, according to Lemma~\ref{A2}, this dependency  can be
removed if we take $i$ large enough.



%
%
%
\bigskip

For obtaining (E$_3$) and (E$_4$) we consider the dynamical system
$\hat{F}=(F\times F)^T:\Delta\times\Delta\circlearrowleft$. It
follows from the definition of the sequence $\{T_n\}$ in \eqref{tt}
that
\begin{equation}\label{tt2}
\hat{F}^n=\left(F\times F\right)^{T_n}, \quad\text{for all $n\ge1$}.
\end{equation}
Let $\hat{\xi}_1$ denote the partition into rectangles $\hat{\Gamma}$ of  $\Delta\times\Delta$
on which  $T$ is constant and  $\hat{F}^n$ maps
$\hat{\Gamma}$ bijectively to $\Delta_0\times\Delta_0$. Next we define inductively partitions $\hat\xi_2,\hat\xi_3,\dots$ of
$\Delta\times\Delta$ by
$
\hat{\xi}_n:=\hat{F}^{-(n-1)}\hat{\xi_1},
$
for $n\ge2$.
Each $\hat{\xi}_n$ is the partition into subsets $\hat{\Gamma}$ of $\Delta\times\Delta$   on which $T_n$ is constant and
$\hat{F}$ maps $\hat{\Gamma}$ bijectively to $\Delta_0\times\Delta_0$. 
We consider the reference measure $m\times m$  for the dynamical system $\hat{F}$
and $J\hat{F}$ the Jacobian of $\hat{F}$ with respect to $m\times m$. We define a separation time
 $\hat{s}:(\Delta\times\Delta)\times (\Delta\times\Delta)\to\NN_0$ for
$\hat{F}$ in the following way: given $w,z\in\Delta\times\Delta$, take
 $$
 \hat{s}(w,z)=\min\big\{n\ge0\colon \text{$\hat{F}^nw$ and $\hat{F}^nz$ lie in distinct elements of $\hat{\xi}_1$}\big \}.
 $$
Denoting
$$
 \Phi=\frac{dP}{d (m\times m)},
 $$
we have $\Phi(x,x')=\varphi(x)\varphi'(x')$. With no loss of
generality we assume from here on  that $\varphi(x)>0$ and
$\varphi(y)>0$. The next two results are proved in \cite[Sublemma
3]{Y2}.

\cle \label{reg J} Let $C_1>0$ be as in Lemma~\ref{l.cfbeta}.
 Given $w,z\in \Delta\times\Delta$ with
$\hat{s}(w,z)\ge n\ge1$
$$\log \frac{J\hat{F}^n(w)}{J\hat{F}^n(z)}\leq
2C_{1}\beta^{\hat{s}(\hat{F}^n(w),\hat{F}^n(z))}.$$
 \fle

\cle \label{reg Phi} Let $C_{\Phi}=C_{\varphi}+C_{\varphi'}$. Given   $w,z\in \Delta\times\Delta$
$$
\log \frac{\Phi(w)}{\Phi(z)}\leq
C_{\Phi}\beta^{\hat{s}(w,z)}.
$$
\fle

The next result gives a distortion control similar to that of Lemma~\ref{A2}.

\cle  \label{A21} There is $C_*=C_*(C_\varphi,C_{\varphi'})>0$ such
that for any  $i\ge 1$  and any $\Gamma\in \hat\xi_i$, we have for
 all $x,y \in \Delta_0\times\Delta_0$ and $Q=\hat F_*^i(P \vert \Gamma)$
$$\left|{\frac{dQ}{dm}(x)}\bigg/{\frac{dQ}{dm}(y)}\right|\leq C_*.$$
 \fle \dem Let $x_0, y_0\in \Gamma$ be such that
$\hat F^i(x_0)=x$ and $\hat F^i(y_0)=y$. Recall that  $\hat{s}(x_0,y_0)\geq i$.
Using Lemma~\ref{reg Phi} and Lemma~\ref{reg J} we obtain
\begin{equation*}
\left|{\frac{dQ}{dm}(x)}\bigg/{\frac{dQ}{dm}(y)}\right|
=\left|\frac{\Phi(x_0)}{J\hat F^i(x_0)}\cdot\frac{J\hat F^i(y_0)}{\Phi(y_0)}\right|
\le  \frac{\Phi(x_0)}{\Phi(y_0)}\cdot \left|\frac{J\hat F^i(y_0)}{J\hat F^i(x_0)}\right|
\leq\exp(C_\Phi+2C_{ 1}).
\end{equation*}
We just have to take $C_*=\exp(C_\varphi+C_{\varphi'}+2C_{1})$. \cqd

 Now we are going to define a sequence of densities
$\hat{\Phi}_0\geq\hat{\Phi}_1\geq\hat{\Phi}_2\geq\cdots$ in $\Delta\times \Delta$ with the property
that for all $i\geq 0$ and all $\hat{\Gamma}\in\hat{\xi}_i$
\begin{equation}\label{dens arrumada}
\pi_*\hat{F}^i_*((\hat{\Phi}_{i-1}-\hat{\Phi}_i)((m\times
m)\vert\hat{\Gamma}))=\pi'_*\hat{F}^i_*((\hat{\Phi}_{i-1}-\hat{\Phi}_i)((m\times m)\vert\hat{\Gamma}))
\end{equation}
Let $\varepsilon=\varepsilon(F)>0$ be a small number to be determined later
(see Lemma~\ref{e1} below). Let $i_1=i_1(\Phi)$ be such that
\begin{equation}\label{equium}
    C_{\Phi}\beta^{i_1}<C_{\hat{F}}.
\end{equation}
 For $i<i_1$, we take
$\hat{\Phi}\equiv\Phi$. For $i\ge i_1$, let
\begin{equation}\label{phich}
\hat{\Phi}_i(z)=\left[\frac{\hat{\Phi}_{i-1}(z)}{J\hat{F}^i(z)}-\varepsilon\min_{w\in\hat{\xi}_i(z)}
\frac{\hat{\Phi}_{i-1}(w)}{J\hat{F}^i(w)}\right]J\hat{F}^i(z),
\end{equation}
where $\hat{\xi}_i(z)$ is the element of $\hat{\Gamma}$ of $\hat{\xi_i}$ that contains $z$. One can easily see
that the sequence
$\{\hat{\Phi}_i\}$ satisfies condition  \eqref{dens arrumada}. 
The next result is proved in \cite[Lemma 3]{Y2}. As observed in
\cite[page 166]{Y2}, $\vare$ depends only on $\beta$.

\cle \label{e1} If $\varepsilon>0$  is sufficiently small, then
there is $0<\varepsilon_1<1$ (not depending on $\Phi$) such that $
\hat{\Phi}_i\leq (1-\varepsilon_1)\hat{\Phi}_{i-1} $ for all $i\geq
i_1$. \fle

\subsubsection{Proof of (E$_3$)} Let $\varepsilon_1>0$ be as in Lemma~\ref{e1}. 
Let $\Phi_0,\Phi_1,\Phi_2,...$ be defined in the following way:
given $n\ge 0$ and  $z\in\Delta\times\Delta$, let
\begin{equation}\label{phin}
\Phi_{n}(z)=\hat{\Phi}_i(z) \quad\textrm{ for }\quad T_i(z)\leq n<T_{i+1}(z).
\end{equation}
 We claim that
 \begin{equation}\label{seikla}
 \mdlo{F_*^n\lambda-F_*^n\lambda'}\leq 2\int\Phi_n d(m\times
m)\quad\text{for all $n\ge 1$}.\end{equation}
 Actually, taking $\displaystyle\Phi=\Phi_n+\sum_{k=1}^n(\Phi_{k-1}-\Phi_k)$ we have
$$\begin{array}{lll}\mdlo{F_*^n\lambda-F_*^n\lambda'}&=&\mdlo{\pi_*(F\times
F)_*^n(\Phi(m\times m))-\pi'_*(F\times F)_*^n(\Phi(m\times
m))}\\
&\leq&\mdlo{\pi_*(F\times F)_*^n(\Phi_n(m\times m))-\pi'_*(F\times F)_*^n(\Phi_n(m\times m))} \\&&+
\displaystyle\sum_{k=1}^n\mdlo{(\pi-\pi')_*[(F\times F)_*^n((\Phi_{k-1}-\Phi_k)(m\times m))]} .\end{array}$$
For the first term in the last sum we have
 $$
 \mdlo{\pi_*(F\times F)_*^n(\Phi(m\times m))-\pi'_*(F\times F)_*^n(\Phi(m\times
m))}\leq 2\int\Phi_nd(m\times m).
 $$
Let us see that all the other terms vanish. Define
 $
 A_{k,i}=\{z\in\Delta\times\Delta: k=T_i(z)\}
 $ and $A_{k}=\cup
A_{k,i}$.
 Each $A_{k,i}$ is a union of elements of
$\Gamma\in\hat{\xi}_i$ and  $A_{k,i}\neq A_{k,j}$ for $i\neq j$. By \eqref{phin} we have
$\Phi_{k-1}-\Phi_k=\hat{\Phi}_{i-1}-\hat{\Phi}_i$ on $\Gamma\in\hat{\xi}_i\vert A_{k,i}$, and
$\Phi_k=\Phi_{k-1}$ on
$\Delta\times\Delta-A_k$. For $k\ge1$
\begin{eqnarray*}
\lefteqn{\hspace{-1cm}\pi_*(F\times F)_*^n((\Phi_{k-1}-\Phi_k)(m\times m))}\\
&=&\sum_i\sum_{\Gamma\subset A_{k,i}}F^{n-k}_*\pi_*(F\times F)_*^{T_i}((\hat{\Phi}_{i-1}-\hat{\Phi}_i)(m\times
m)\vert\Gamma)\\
&=&\sum_i\sum_{\Gamma\subset A_{k,i}}F^{n-k}_*\pi'_*(F\times
F)_*^{T_i}((\hat{\Phi}_{i-1}-\hat{\Phi}_i)(m\times m)\vert\Gamma), \quad \quad \text{by \eqref{dens arrumada}}\\
&=&\pi'_*(F\times F)_*^n((\Phi_{k-1}-\Phi_k)(m\times m)).
\end{eqnarray*}
This completes the proof of \eqref{seikla}. To finish (E$_3$) we
write
$$\int \Phi_n d(m\times m)= \int_{\{T_{i_1}>n\}}\Phi_n d(m\times m)+
\sum_{i=i_1}^{\infty}\int_{\{T_i\leq n< T_{i+1}\}}\Phi_n d(m\times m).$$
Observe that
$$\int_{\{T_{i_1}>n\}}\Phi_n d(m\times m)=\int_{\{T_{i_1}>n\}}\Phi d(m\times m)=P\{T_{i_1}>n\},$$
while for $i\geq i_i$,
\begin{eqnarray*}
\lefteqn{\hspace{-3cm}\int_{\{T_i\leq n< T_{i+1}\}}\Phi_n} d(m\times m)&=&\int_{\{T_i\leq n< T_{i+1}\}}\hat{\Phi}_i d(m\times m)\\[0.3cm]
&\leq&\int_{\{T_i\leq n<
T_{i+1}\}}(1-\varepsilon_1)^{i-i_1+1}\Phi d(m\times m)\\[0.3cm]
&=&(1-\varepsilon_1)^{i-i_1+1}P\{T_i\leq n< T_{i+1}\}.
\end{eqnarray*}
Hence
\begin{align}
\mdlo{F_*^n\lambda-F_*^n\lambda'} &\leq
2P\{T_{i_1}>n\}+2\sum_{i=i_1}^{\infty}(1-\varepsilon_1)^{i-i_1+1} P\{T_i\leq n< T_{i+1}\}\nonumber\\
&\leq
2P\{T_{i_1}>n\}+2(1-\varepsilon_1)^{-i_1+1}\sum_{i=i_1}^{\infty}(1-\varepsilon_1)^{i} P\{T_i\leq n< T_{i+1}\}. \label{equild}
\end{align}
We may write
\begin{align*}
P\{T_{i_1}>n\}&= P\{T>n\} +(1-\varepsilon_1)^{-i_1+1}\sum_{i=1}^{i_1-1}(1-\varepsilon_1)^{i_1-1} P\{T_i\leq n< T_{i+1}\}\\
 &\le P\{T>n\}+ (1-\varepsilon_1)^{-i_1+1}\sum_{i=1}^{i_1-1}(1-\varepsilon_1)^{i} P\{T_i\leq n< T_{i+1}\},
\end{align*}
 which together with \eqref{equild} yields
 $$
\mdlo{F_*^n\lambda-F_*^n\lambda'}\le 2P\{T>n\}+ K_1 \sum_{i=1}^{\infty}(1-\varepsilon_1)^i P\{T_i\leq n < T_{i+1}\},
 $$
 with $K_1$ depending only on $\vare_1$ and $i_1$. From \eqref{equium} and Lemma~\ref{e1} one easily obtains the desired dependence of
 $K_1$ on $\varphi$ and $\varphi'$.

\subsubsection{Proof of (E$_4$)}
This estimate is obviously true for $i=0$. Take an arbitrary  $i\ge
1$ and $\Gamma\in \hat\xi_i$. Recall that $\hat{F}^i$ maps
$\hat{\Gamma}$ bijectively to $\Delta_0\times\Delta_0$. Letting
$Q=F^{i}_*(P\vert \Gamma)$ and observing that from \eqref{tt} and
\eqref{tt2} we have
  $T_{i+1}-T_i=T\circ \hat F^i$, we may write
$$P\{T_{i+1}-T_i>n\mid{\Gamma}\}=\frac{Q\{T>n\}}{Q(\Delta_0\times\Delta_0)}.$$
Using Lemma~\ref{A21},
$$P\{T_{i+1}-T_i>n\mid{\Gamma}\}\leq C_*^2\frac{(m\times m)\{T>n\}}{(m\times m)(\Delta_0\times\Delta_0)}.$$
From this last inequality one easily obtains (E$_4$) with
$K_2=C_*^2/(m\times m)(\Delta_0\times\Delta_0)$.

\end{document}